\documentclass[11pt]{article}
\usepackage{amsmath,amscd,amsthm,amssymb,latexsym,verbatim}

\usepackage{graphicx,epsf}

-\marginparwidth \oddsidemargin -4mm \evensidemargin \oddsidemargin

\newtheorem{theorem}{Theorem}[section]

\newtheorem{proposition}[theorem]{Proposition}
\newtheorem{lemma}[theorem]{Lemma}
\newtheorem{corollary}[theorem]{Corollary}
\theoremstyle{definition}

\theoremstyle{remark}

\newcommand{\commentout}[1]{}

\newcounter{smalllist}

\allowdisplaybreaks
\numberwithin{equation}{section}

\newcommand{\abs}[1]{\lvert#1\rvert}

\newcommand{\no}{\nonumber}
\newcommand{\lb}{\label}

\newcommand{\supp}{\text{\rm{supp}}}

\newcommand{\beq}{\begin{equation}}
\newcommand{\eeq}{\end{equation}}

\newcommand{\bal}{\begin{align}}
\newcommand{\eal}{\end{align}}
\newcommand{\bals}{\begin{align*}}
\newcommand{\eals}{\end{align*}}

\newcommand{\bbR}{{\mathbb{R}}}
\newcommand{\RR}{{\mathbb{R}}}

\newcommand{\calS}{{\mathcal S}}

\newcommand{\eps}{\varepsilon}

\newcommand{\tht}{\theta}

\newcommand{\farc}{\frac}
\newcommand{\Rm}{{\mathbb R}}
\newcommand{\e}{{\varepsilon}}
\newcommand{\un}{{\mathbf{1}}}

\def\di{\displaystyle}
\def\supp{\mathrm{supp}}
\newcommand{\lambdam}{{\lambda'}}
\newcommand{\blambda}{{\lambda}}
\newcommand{\galpha}{{\alpha}}

\begin{document}
\title{Existence and Non-existence of Fisher-KPP Transition Fronts}
\begin{author}{James Nolen\thanks{Department of Mathematics, Duke University, 
Durham, NC 27708;
nolen@math.duke.edu} \and 
Jean-Michel Roquejoffre\thanks{Institut de Math\'ematiques, (UMR CNRS 5219), 
Universit\'e Paul Sabatier, 118 route de Narbonne, 31062 Toulouse cedex, France; 
roque@mip.ups-tlse.fr} \and Lenya Ryzhik\thanks{Department of Mathematics,  
Stanford University, Stanford CA 94305;
ryzhik@math.stanford.edu} \and 
Andrej Zlato\v s\thanks{Department of Mathematics,
University of Wisconsin,
480 Lincoln Drive,
Madison, WI 53706, USA; \hbox{andrej@math.wisc.edu}}}
\end{author}
\maketitle

\begin{abstract}
We consider Fisher-KPP-type reaction-diffusion equations with spatially 
inhomogeneous
reaction rates. We show that   a sufficiently strong localized inhomogeneity may 
prevent existence of 
transition-front-type global in time solutions while creating a global in time 
bump-like solution.  This is the first example of a medium in which no reaction-diffusion transition front exists. 
A weaker localized inhomogeneity leads to existence of transition fronts but 
only in a finite range
of speeds. These results are in contrast with  both Fisher-KPP reactions in
homogeneous media 
as well as  ignition-type reactions in inhomogeneous media.
\end{abstract}
\section{Introduction and main results}

\subsection*{Fisher-KPP traveling fronts in homogeneous media}

Traveling front solutions of the reaction-diffusion equation
\begin{equation}\label{apr2801}
u_t=u_{xx}+f(u)
\end{equation}
are used to model phenomena
in a range of applications from biology to social sciences, and have been studied extensively since the pioneering papers of Fisher~\cite{Fisher} and Kolmogorov-Petrovskii-Piskunov~\cite{KPP}.  The Lipschitz nonlinearity  $f$ is said to be of {\it KPP-type} if 
\begin{equation}
\label{apr2802}
f(0)=f(1)=0 \qquad\text{and}\qquad 0<f(u)\le f'(0)u \text{ \,\,for } u\in(0,1),
\end{equation}
and one considers solutions $0< u(t,x)< 1$.
A {\it traveling
front} is a solution of (\ref{apr2801}) of the form $u(t,x)=\phi_c(x-ct)$, with the 
function $\phi_c(\xi)$ satisfying
\begin{equation}\label{apr2804}
\phi_c''+c\phi_c'+f(\phi_c)=0,~~\phi_c(-\infty)=1,~~\phi_c(+\infty)=0.
\end{equation}
Here $c$ is the {\it speed} of the front and traveling fronts exist precisely when $c\ge c_*\equiv 2\sqrt{f'(0)}$. For the sake of convenience we will
assume that $f'(0)=1$, which can be achieved
by a simple rescaling of space or time. 

The traveling front profile $\phi_c(\xi)$  satisfies $\phi_c(\xi)\sim e^{-r(c)\xi}$ as $\xi\to +\infty$. The decay rate $r(c)$ can be obtained from the linearized  problem $v_t=v_{xx}+v$, and is given by
\begin{equation}\label{apr2806}
r(c)=\frac{c-\sqrt{c^2-4}}{2}.
\end{equation}
It is the root of both $r^2-cr+1=0$ and $r^2+r\sqrt{c^2-4}-1=0$
and for $c\gg 1$ we have 
$r(c)=c^{-1}+O(c^{-3})$,
whence $\lim_{c\to+\infty}cr(c)=1$.

\subsection*{Fisher-KPP transition fronts in inhomogeneous media and  bump-like solutions}

In this paper we consider the {\it inhomogeneous} 
reaction-diffusion equation
\beq \lb{1.1}
u_t=u_{xx} + f(x,u)
\eeq
with $x\in\bbR$ and a KPP reaction $f$. That is, we assume that $f$ is Lipschitz, $f_u(x,0)$ exists, 
\begin{equation}\label{apr2812}
f(x,0)=f(x,1)=0, \qquad\text{and} \qquad  0< f(x,u)\le f_u(x,0)u \text{ \,\,for $(x,u)\in\bbR\times(0,1)$.}
\end{equation} 
We let $a(x)\equiv f_u(x,0)>0$ and assume that for some $C,\delta>0$ we have
\beq \lb{10.1}
f(x,u)\ge a(x)u - Cu^{1+\delta} \text{ \,\,for $(x,u)\in\bbR\times(0,1)$}.
\eeq
Finally, we will assume here
\begin{equation}\label{apr3008}
\hbox{$0<a_-\le a(x)\le a_+<+\infty$ \,\,for  $x\in\Rm$}
\end{equation}
and
\begin{equation}\label{apr3004}
\hbox{$\lim_{|x|\to\infty} a(x)= 1$.}
\end{equation}
That is, we will consider media which are localized perturbations of the homogeneous case.

In this case traveling fronts with a constant-in-time profile cannot exist in general, and one instead considers {\it transition fronts}, a generalization of traveling fronts introduced in~\cite{BH-Brezis,Matano,Shen1}. In the present context, 
a global in time solution of (\ref{1.1}) is said to be a transition front if 
\begin{equation}\label{apr2814}
\lim_{x\to-\infty}u(t,x)=1 \qquad\text{and}\qquad \lim_{x\to+\infty}u(t,x)=0
\end{equation}
for any $t\in\bbR$, and for any $\eps>0$ there exists $L_\eps<+\infty$ such that for any $t\in\bbR$ 
we have 
\begin{equation}\label{apr2816}
\hbox{diam}\left\{x\in\Rm \,|\,\eps\le u(t,x)\le 1-\eps\right\}<L_\eps.
\end{equation}
That is, a transition front is a global in time solution connecting $u=0$ and $u=1$ at any time $t$, which also has a uniformly 
bounded in time width of the transition region between $\eps$ and $1-\eps$.

Existence of transition fronts has been previously established for a class of 
time-dependent spatially homogeneous {\it bistable} nonlinearities
in~\cite{Shen1}, and for spatially inhomogeneous {\it ignition}  nonlinearities 
in~\cite{MRS,NolR,Zlatos}. The results in these papers, while non-trivial, are similar in spirit to the situation for such  
nonlinearities in homogenous media: there exists a unique
(up to a time shift) transition front, and it is asymptotically stable for 
 the Cauchy problem.  In the present paper we will demonstrate that the situation can be very different for  
KPP-type nonlinearities, even in the case of a spatially localized inhomogeneities.


Before we do so, let us define another type of a solution of \eqref{1.1}.
We say that a global in time solution $0<u(t,x)<1$ of (\ref{1.1}) is {\it bump-like} 
if  $u(t,\cdot)\in L^1(\Rm)$ for all $t\in\Rm$. We will show that bump-like solutions can exist for inhomogeneous KPP-type nonlinearities. What makes such solutions special is that they do not exist in many previously studied settings, as can be seen from the following proposition. 

\begin{proposition}
\label{p1.1}
Assume that either  $f(x,u)\ge 0$ is an ignition reaction (i.e.,  $f(x,u)=0$ if  $u\in [0,\theta(x)]\cup\{1\}$, with $\tht\equiv \inf_{x\in\bbR}\theta(x)>0$; see \cite{MRS,NolR, Zlatos}) or $f(x,u)=f(u)$ is a spatially homogeneous KPP reaction satisfying \eqref{apr2802} and
\begin{equation}\label{apr2902x}
\hbox{$f(u)\equiv u\,\,$  for $u\in[0,\theta]$}
\end{equation} 
for some $\theta\in(0,1)$.
Then \eqref{1.1} does not admit global in time bump-like solutions.  
\end{proposition}

{\it Remark.}  Hypothesis \eqref{apr2902x} is likely just technical but we make it for the sake of simplicity.
\smallskip

\subsection*{Non-existence of transition fronts for strong KPP inhomogeneities}

Our first main result shows that a localized KPP inhomogeneity can create global in time 
bump-like solutions of (\ref{1.1}) as well as prevent existence of any transition front solutions.  This is the first example of a medium in which no reaction-diffusion transition fronts exist.  Moreover, in the case $a(x)\ge 1$ and $a(x)-1$ compactly supported, Theorems \ref{T.1.1} and \ref{T.1.2} together provide a sharp criterion for the existence of transition fronts. Namely, transition fronts exist when $\blambda<2$ and do not exist when $\blambda>2$, with
 $\blambda\equiv \sup \sigma(\partial_{xx} + a(x))$ the supremum of the spectrum of the operator $L\equiv \partial_{xx} + a(x)$ on $\bbR$.
One can consider these to be the main results of this paper.

Note that \eqref{apr3004} implies that the essential 
spectrum of $L$ is $(-\infty,1]$ and so $\blambda\ge 1$. Hence if $\blambda>1$ 
then $\blambda$ is the principal eigenvalue of $L$ and
\begin{equation}\label{apr2904}
\psi''+a(x)\psi=\blambda\psi
\end{equation}
holds for the positive eigenfunction~$0<\psi\in L^2(\Rm)$ satisfying also $\|\psi\|_\infty=1$.
We note that  $\psi(x)$ decays exponentially as $x\to \pm\infty$ due to \eqref{apr3004}.

\begin{theorem} \lb{T.1.1}
Assume that $f(x,u)$ is a KPP reaction satisfying 
\eqref{apr2812}--\eqref{apr3004} with $a_-= 1$.
If $\blambda>2$, then any global in time 
solution of \eqref{1.1}  such that $0<u(t,x)<1$ satisfies (with $C_{c}>0$)
\beq \lb{10.2}
u(t,x)\le C_{c} \,e^{-|x|+ct}
\eeq
for any $c< \blambda/\sqrt{\blambda-1}$  and all $(t,x)\in\bbR^-\times\bbR$.  In particular, no transition front exists.   

Moreover, bump-like solutions do exist, and if there is $\tht>0$ such that
\begin{equation}\label{apr2902}
f(x,u)\equiv a(x)u  \text{ \,\,for all $(x,u)\in\bbR\times[0,\theta]$},
\end{equation}
then there is a unique (up to a time-shift) global in time solution $0<u(t,x)<1$. This solution satisfies $u(t,x)= e^{\blambda t} \psi(x)$ for $t\ll -1$.
\end{theorem}

\subsection*{Existence and non-existence of  transition fronts for weak KPP inhomogeneities}

We next show that transition fronts do exist when $\blambda<2$, albeit in a bounded range of speeds. If $u$ is a transition front, let $X(t)$ be the rightmost point $x$ such that $u(t,x)=1/2$. If
$$
\lim_{t-s\to+\infty}\frac{X(t)-X(s)}{t-s}=c,
$$
then we say that $u$ has {\it global mean speed} (or simply {\it speed}) $c$.
Recall that in the homogeneous KPP case with $f'(0)=1$, traveling fronts exist for all speeds $c\ge 2$.  

\begin{theorem} \lb{T.1.2}
Assume that $f(x,u)$ is a KPP reaction satisfying 
\eqref{apr2812}--\eqref{apr3004} and $a(x)-1$ is compactly supported. If $\blambda\in(1,2)$, then  
for each $c\in(2,\blambda/\sqrt{\blambda -1})$ equation \eqref{1.1} admits 
a transition front solution  with global mean
speed $c$. Moreover, bump-like solutions also exist.
\end{theorem}

{\it Remarks.}  1. In fact, the constructed fronts will satisfy  $\sup_{t\in\bbR}|X(t)-ct|<\infty$.
\smallskip


2.  Fisher-KPP  equations in homogeneous media also
admit global in time solutions that are mixtures of traveling
fronts moving with different speeds, constructed in~\cite{HN1,HN2}. 
Such global in time mixtures of transition fronts constructed in Theorem~\ref{T.1.2} 
also exist, but this problem will be considered elsewhere in order to keep this paper concise. Existence of transition fronts with the critical 
speeds $c_*=2$ and 
$c^*\equiv\blambda/{\sqrt{\blambda-1}}$ is a delicate issue and will also be left for a later work.
\smallskip


Finally, we show that the upper limit $\blambda/\sqrt{\blambda -1}$ on the front speed in Theorem \ref{T.1.2} is not due to our techniques being inadequate. Indeed, we will prove non-existence of fronts with speeds $c>\blambda/\sqrt{\blambda -1}$, at least under additional, admittedly somewhat strong, conditions on $f$.  

\begin{theorem} \lb{T.1.3}
Assume that $f(x,u)=a(x)f(u)$ where $a$ is even, satisfies  \eqref{apr3008} with $a_-= 1$, and $a(x)-1$ is compactly supported, and $f$ is such that \eqref{apr2802} and \eqref{apr2902x} hold for some  $\tht\in(0,1)$.
In addition assume that \eqref{apr2904} has a unique eigenvalue $\blambda > 1$. Then there are no transition fronts  with global mean speeds  $c> \blambda/{\sqrt{\blambda-1}}$.
\end{theorem}



Let us indicate here the origin of the threshold $\blambda/{\sqrt{\blambda-1}}$ for speeds of transition fronts.  In the homogeneous case $f(x,u)=f(u)$ with $f(u)=u$ for $u\le\tht$, the traveling front with speed $c\ge 2$ satisfies $u(t,x)=e^{-r(c)(x-ct)}$ (up to a time shift) for $x\gg ct$. This means that $u$ increases at such $x$ at the exponential rate $cr(c)$ in $t$.  We have $\lim_{|x|\to\infty} f_u(x,0)= 1$, so it is natural to expect a similar behavior of a transition front $u$ (with speed $c$) at large $x$.  On the other hand, any non-negative non-trivial solution of \eqref{1.1} majorizes a multiple of $e^{\lambda_Mt}\psi_M(x)$ for $t\ll -1$, with $\lambda_M$ and $\psi_M$ the principal eigenvalue and eigenfunction of $\partial_{xx} + a(x)$ on $[-M,M]$ with Dirichlet boundary conditions (extended by 0 outside $[-M,M]$). So $u$ has to increase at least at the rate $\lambda_M$, and since $\lim_{M\to\infty}\lambda_M=\blambda$, it follows that one needs $cr(c)\ge \lambda$ in order to expect existence of a transition front with speed $c$. Using \eqref{apr2806}, this translates into $c\le \blambda/{\sqrt{\blambda-1}}$.

In the rest of the paper we prove Proposition \ref{p1.1} and Theorems  \ref{T.1.1},  \ref{T.1.2}, \ref{T.1.3}  (in Sections 2, 3, 4, and 5--7, respectively).
\medskip

{\bf Acknowledgment.} JN was supported by NSF grant DMS-1007572, 
JMR by ANR grant 'PREFERED', LR by NSF grant DMS-0908507, and AZ by
NSF grant DMS-0901363 and an Alfred P. Sloan Research Fellowship.

\section{Nonexistence of bump-like solutions for ignition reactions and homogeneous 
KPP reactions: The proof of Proposition \ref{p1.1}.}

Assume, towards contradiction, that there exists a bump-like solution.  We note that parabolic regularity and $f$ Lipschitz then yield for each $t\in\bbR$,
\[
u,u_x \to 0 \qquad \text{as $|x|\to \infty$}.
\]
This will guarantee that differentiations in $t$ of integrals over $\bbR$ and  integration by parts below are valid.
Let us define
\[
I(t) \equiv \int_\RR u(t,x)\, dx \qquad \text{and} \qquad J(t) \equiv \frac 12 \int_\RR u(t,x)^2\, dx.
\]
Integration of \eqref{1.1} and of \eqref{1.1} multiplied by $u$ over $x\in\RR$  yields
\[
I'(t) = \int_\RR f(x,u)\, dx \ge 0 \qquad \text{and} \qquad J'(t) = \int_\RR f(x,u)u\, dx - \int_\bbR |u_x|^2 \,dx \le I'(t) - \int_\bbR |u_x|^2 \,dx.
\]
So $\lim_{t\to-\infty} I(t)=C\ge 0$ and then $\lim_{t\to-\infty} \int_\bbR |u_x|^2 \,dx = 0$. Parabolic regularity again gives  
\[
u,u_x \to 0 \qquad \text{as $t\to -\infty$, uniformly in $x$}.
\]
Thus $u(x,t)\le \theta$ for all $t<t_0$ and all $x\in\bbR$.  Then $u$ in the ignition case ($v(t,x)\equiv e^{-t}u(t,x)$ in the KPP case) solves the heat equation for $t\le t_0$. 
 Since $u\ge 0$ ($v\ge 0$) and it is $L^1$ in $x$, it follows that $u=0$  ($v=0$), a contradiction.

 \section{The case $\blambda>2$: The proof of Theorem~\ref{T.1.1}}\label{sec:2}

We obviously only need to consider $c\in (2, \blambda/\sqrt{\blambda-1})$, so let us assume this. We will first assume, for the sake of simplicity, that $a(x)-1$ is compactly supported and \eqref{apr2902} holds. At the end of this section we will show how to accommodate the proof to the general case. 

Let us shift the origin by a large enough $M$  so that in the
shifted coordinate frame $a(x)\equiv 1$ for $x\notin[0,2M]$, and  
the principal eigenvalue $\lambda_M$ of $\partial_{xx}+a(x)$ on $(0,2M)$ with Dirichlet boundary conditions  satisfies  $\lambda_M >2$.  This is possible since
\[
\lim_{M\to+\infty}\lambda_M = \blambda.
\]
 We let $\psi_M$ be the corresponding $L^\infty$-normalized principal eigenfunction, that is,  $\lVert \psi_M \rVert_\infty = 1$ and
\begin{equation}\label{apr2901}
\psi_M''+a(x)\psi_M=\lambda_M\psi_M,~~\psi_M>0
\hbox{ on $(0,2M)$},~~\psi_M(0)=\psi_M(2M)=0.
\end{equation}

It is easy to show that any entire solution $u(t,x)$ of~(\ref{1.1})
such that $0<u(t,x)<1$ satisfies
$\lim_{t\to-\infty} u(t,x) =0$ and 
$\lim_{t\to +\infty} u(t,x) =1$ for any $x\in\bbR$, so after a possible translation of $u$ forward in time by some $t_0$, we can assume 
\beq \lb{1.1a}
\sup_{t\le 0} u(t,M) <  \tht \psi_M(M) \le  \tht.
\eeq
In that case \eqref{10.2} for this translated $u$  yields $u(t,x)\le C e^{-|x-M|+c(t-t_0)}$ when $t< t_0$ for the original $u$, but then the result follows for a larger $C$ from the fact that $C e^{-|x-M|+(1+\|a\|_\infty)(t-t_0)}$ is a supersolution of \eqref{1.1} on $(-t_0,0)\times\bbR$.

\subsection*{Non-existence of transition fronts}

Assume that $u$ is a global in time solution of \eqref{1.1}. Non-existence of transition fronts obviously follows from \eqref{10.2}.  
The following lemma is the main step in the proof of (\ref{10.2}). 

\begin{lemma}\label{lem1apr30}
For any $c,c'\in(2,\lambda_M/\sqrt{\lambda_M-1})$ with $c<c'$,  there is $C_0>0$ (depending only on $a$, $\tht$, $c$, $c'$) and $\tau_0>0$ (depending also on $u(0,M)$) such that 
\begin{equation}\label{apr3020}
u(t,x)\le C_0u(0,M)e^{x+ct}
\end{equation}
holds  for all $t\le -1$ and $x\in[0,c'(-t-1)]$, as well as  for all $t\le -\tau_0$ and $x\ge 0$.
\end{lemma}

{\it Remark.} This is a one-sided estimate but by symmetry of the arguments in its proof, the same estimate holds for $u(-t,2M-x)$. 
\smallskip

Let us show how this implies (\ref{10.2}), despite the fact that (\ref{apr3020})
seemingly goes in two wrong directions. First, the estimate holds for $x\ge 0$ but the exponential on the right 
side grows as $x\to+\infty$. Second, this exponential
is moving to the left as time progresses in the positive direction, while we are estimating $u$ to the right of $x=0$.  
The point of \eqref{apr3020} is that the speed $c$ at which  the exponential moves is larger than 2, the latter being the minimal speed of fronts when $a(x)=1$ everywhere.  Thus, when  looking at large negative times, this gives us a much smaller than expected upper bound on $u$ at $|x|\le c|t|$.  Using this bound and then going forward in time towards $t=0$, we will find that $u$ cannot become $O(1)$ at $(0,M)$.


Given $c\in (2, \blambda/\sqrt{\blambda-1})$, pick $M$ such that $c< \lambda_M/\sqrt{\lambda_M-1})$ and then $c'>c$ as in Lemma \ref{lem1apr30}. Let $\tau_1\equiv 1+2M /c' $ (so $\tau_1$ depends on $a,\tht,c$ but not on $u$). By the first claim of Lemma \ref{lem1apr30} we have
\begin{equation}\label{apr3028}
u(t,2M)\le C_0u(0,M)e^{2M+ct}
\end{equation}
for all $t\le -\tau_1$ because then $2M\le c'(-t-1)$.

Next, for any $t_0\le -\tau_0$, we let
\[
v_{t_0}(t,x) \equiv C_0 u(0,M)  e^{x+ct_0 + 2(t-t_0)} + C_0 u(0,M) e^{4M-x+ct}.
\]
Then $v_{t_0}$ is a super-solution for  \eqref{1.1} on $(t_0,\infty)\times(2M,\infty)$ 
since $a(x)\equiv 1$
for $x> 2M$. Moreover, the second claim of Lemma \ref{lem1apr30} and $t_0\le -\tau_0$  imply that
at the ``initial time'' $t_0$ we have
\[
u(t_0,x)\le C_0 u(0,M)  e^{x +c t_0 }\le v_{t_0}(t_0,x)
\]
for all $x>2M$. 
Since $c>2$, it follows from (\ref{apr3028}) 
that $u(t,2M)\le v_{t_0}(t,2M)$
 for all $t\in(t_0,-\tau_1)$.
Since the super-solution $v_{t_0}$ is above $u$ initially (at $t=t_0$) on all 
of $(2M,\infty)$ and at $x=2M$ for all $t\in(t_0,-\tau_1)$, the maximum principle 
yields 
\begin{equation}\label{apr3022}
u(t,x) \le v_{t_0}(t,x) 
\end{equation}
for all $t\in[t_0,-\tau_1]$ and $x\ge 2M$. Since $c>2$, taking  $t_0\to -\infty$ in (\ref{apr3022}) gives
\begin{equation}\label{apr3024bis}
u(t,x)\le C_0 u(0,M)  e^{4M-x+ct},
\end{equation}
for $t\le -\tau_1$ and $x\ge 2M$. 
Note that unlike our starting point (\ref{apr3020}),
the estimate (\ref{apr3024bis}) actually goes in the right
direction, since the exponential is decaying as $x\to+\infty$.

An identical argument gives $u(t,x)\le C_0 u(0,M)  e^{2M+x+ct}$ for $t\le -\tau_1$ and $x\le 0$, so 
\begin{equation}\label{apr3024-1}
u(t,x)\le C_0e^{2M} u(0,M)  e^{-|x|+ct}
\end{equation}
for $t\le -\tau_1$ and $x\in \bbR\setminus (0,2M)$. 
Harnack inequality extends this bound to all $t\le -\tau_1-1$ and $x\in\bbR$,
with some $C_1$ (depending only on $a$ and $\tht$) in place of $C_0e^{2M}$:
\begin{equation}\label{apr3024}
u(t,x)\le C_1  u(0,M)  e^{-|x|+ct}
\end{equation}
for all $t\le-\tau_1 -1$ and $x\in\bbR$.
Finally, it follows from (\ref{apr3024}) that
\[
u(t,x)\le C_1 u(0,M)  e^{-|x|+c(-\tau_1-1)} e^{(1+\|a\|_\infty)(t-(-\tau_1-1))}
\]
for $t\ge -\tau_1-1$ because the right-hand side is a super-solution of \eqref{1.1}. 
Since $\tau_1$ only depends on $a,\tht,c$ (once $M,c'$ are fixed) and not on $u$, and since $a_1\ge 1$,
it follows that
\beq\lb{1.7}
u(t,x)\le C_2 u(0,M) e^{-|x|+ct}
\eeq
for all $t\le 0$ and $x\in \bbR$, with $C_2$ depending only on $a,\tht,c$. This is (\ref{10.2}), proving non-existence of transition fronts when $\blambda>2$ under the additional assumptions of $a(x)-1$ compactly supported and \eqref{apr2902} 
(except for the proof of Lemma~\ref{lem1apr30} below). 

\subsection*{Bump-like solutions and uniqueness of a global in time solution}

Existence of a bump-like solution is immediate from \eqref{apr2902}.  Indeed, it is obtained by continuing the solution of \eqref{1.1}, given by  $u(t,x)= e^{\blambda t} \psi(x)$ for $t\ll -1$, to all $t\in\bbR$.

In order to prove the uniqueness claim, we note that the same argument as above, with $u(0,M)$ replaced by $u(s,M)$ and $t\le s\le 0$, gives (with the same $C_2$)
\beq \lb{1.4}
u(t,x)\le C_2 u(s,M)  e^{-|x|-2(s-t)}.
\eeq 
We also have $\|u(t,\cdot)\|_\infty\le\tht$ for all $t\le t_0\equiv - \tfrac 12 \log C_2$.
Therefore, the function $v(t,x)\equiv u(t,x)e^{-2t}$ solves the linear equation 
\begin{equation}\label{may301}
v_t = v_{xx} +(a(x)-2)v
\end{equation}
on $(-\infty, t_0)\times\bbR$. 
It can obviously be extended to an entire solution of (\ref{may301}) by 
propagating it forward in time. 
Taking $t=s$ in \eqref{1.4} gives 
$v(t,x)\le C_2 v(t,M) $ for $(t,x)\in (-\infty,  t_0)\times\bbR$.
Moreover, it is well known that since 
$\blambda$ is an isolated eigenvalue (because 
$\blambda>1$ and the essential spectrum is $(-\infty,1]$), 
the function
$e^{-(\blambda-2) t} v(t,x)$ converges uniformly to 
$\psi(x)$ as $t\to\infty$. 
It follows that 
\beq\lb{1.9}
v(t,x)\le C_3 v(t,M)
\eeq
holds for some $C_3>0$ and all $(t,x)\in\bbR^2$.

We can now apply Proposition 2.5 from~\cite{HP} to (\ref{may301}). 
More precisely, as $a(x)\equiv 1$ outside of a bounded interval, Hypothesis A of 
this proposition is satisfied,
while $\blambda>2$ ensures that Hypothesis H1 of \cite{HP} holds for the 
solution $w(t,x)=e^{(\blambda-2)t} \psi(x)$ 
of (\ref{may301}). Finally, (\ref{1.9}) guarantees that condition 
(2.12) of \cite{HP} holds, too.  It then follows from the
aforementioned proposition that $w(t,x)$ is the unique (up to a time shift) 
global in time solution of (\ref{may301}), 
proving the uniqueness claim in Theorem~\ref{T.1.1}. 

It remains now only to 
prove Lemma~\ref{lem1apr30}
in order to finish the proof of Theorem~\ref{T.1.1}  in the case when $a(x)-1$ is compactly supported and \eqref{apr2902} holds.


\subsection*{The proof of Lemma~\ref{lem1apr30}}

We will prove Lemma~\ref{lem1apr30} 
using the following lemma.

\begin{lemma}\label{L.1.3}
For every $\eps\in(0,1)$ there exists  $C_\eps\ge 1$  (depending also on $a$, $\tht$, and $\lambda_M$)
such that 
\beq \lb{1.2}
u(t,x) \le C_\eps u(0,M)\sqrt{|t|\,}\, e^{\sqrt{\lambda_M  -1}\,x +(\lambda_M   -\eps) t}
\eeq
holds for all $t\le -1$ and $x\in[0,c_\eps (-t-1)]$, with $c_\eps\equiv{(\lambda_M 
-\eps)}/{\sqrt{\lambda_M  -1}} $.
\end{lemma}

Let us first explain how Lemma~\ref{L.1.3} implies 
Lemma~\ref{lem1apr30}.
Pick $\eps>0$ such that 
$c_\eps=c'$.
Then there is $C_0>0$ depending only on $a$, $\tht$, $c$ (via $\eps,\lambda_M,C_\eps$) such that for all $t\le -1$ and $x\in[0,c' (-t-1)]$ we have
\beq \lb{1.3}
u(t,x) \le C_\eps u(0,M)\sqrt{|t|\,}\, e^{\sqrt{\lambda_M  -1}(x +c' t)}  \le C_\eps u(0,M)\sqrt{|t|\,}\, e^{x +c' t}  \le  C_0 u(0,M)  e^{x +c t },
\eeq
the first claim of Lemma~\ref{lem1apr30} .

Next let 
\beq\lb{1.8}
\tau_0 \equiv \frac {|\log(C_0u(0,M)e^{-c})|} {c'-c} +1,
\eeq
so that $C_0 u(0,M)  e^{x +c t}\ge 1$ for $t\le -\tau_0$ and
$x\ge c'(-t-1)$. Since $u(t,x)\le 1$, this means that \eqref{apr3020} also holds for all $t\le -\tau_0$ and $x\ge 0$, the second claim of Lemma~\ref{lem1apr30}.

Thus we are left with the proof of Lemma~\ref{L.1.3}. This, in turn, relies on the following lemma.


\begin{lemma} \lb{L.1.2}
For each $m\in\bbR$ and $\eps>0$ there is $k_\eps>0$ 
such that if $u\in[0,1]$ solves \eqref{1.1} with 
$u(0,x)\ge \gamma \chi_{[l-1,l]}(x)$ for some $\gamma\le \tht/ 2$ and 
$l\in\bbR$, then for $t\ge 0$ and $x\le l+m-2t$,
\[ 
u(t,x) \ge  k_\eps\gamma e^{(1-\eps) t}\int_{l-1}^{l} \frac {e^{-{|x-z|^2}/ {4 
t}}}{\sqrt {4\pi  t}} 
  \, dz.
\]
\end{lemma}

\noindent {\bf Proof.}
The result, with 1 in place of $1-\eps$, clearly holds when $f(x,u)\ge u$ for 
all $x,u$. Since 
$f(x,u)\ge u$ only for $u\le \theta$, we will have to be a little more careful.


It is obviously  sufficient to consider $l=0$. Let $g$ be a concave function on 
$[0,1]$ such that $g (w)=w$ 
for $w\in[0,1/2]$ and $g(1)=0$ and define  $g_\gamma(w)\equiv 2\gamma 
g(w/2\gamma)$ 
(hence $g_\gamma(w)=w$ for $w\in[0,\gamma]$, and $g_\gamma\le f$).
The comparison principle implies that $u(x)\ge w(x)$, where
$w(x)$ solves
\beq \lb{1.5}
w_t=w_{xx}+g_\gamma(w)
\eeq
with initial condition $w(0,x)=\gamma \chi_{[-1,0]}(x)$. It follows from standard 
results on spreading of solutions to 
KPP reaction-diffusion equations (see, for instance, \cite{AW}) that for each 
$\eps>0$ there exists 
$t_\eps\ge {(m+1)}/ {2\sqrt{1-\eps}}$ such that for all  $t\ge t_\eps$ we have 
$w(t,-2\sqrt{1-\eps}(t-t_\eps)-1)\ge \gamma$. 
The time $t_\eps$ is independent of $\gamma$ because  
$w/\gamma$ is independent of $\gamma$.

Note that the function
\[
v(t,x) = e^{-2t_\eps} \gamma e^{(1-\eps)t}\int_{-1}^{0} \frac 
{e^{-|x-z|^2/4 t} }{\sqrt {4\pi  t}}  \, dz
\]
solves $v_t=v_{xx}+(1-\eps)v$, 
so $v$ is a sub-solution of \eqref{1.5} on any domain where $v(t,x)\le \gamma$. 
We have $\| v(t,\cdot)\|_\infty\le e^{-(1-\eps)t_\eps}\gamma \le \gamma$ for 
$t\le t_\eps$, as well as
\[
v(t,-2\sqrt{1-\eps}(t-t_\eps)-1) \le e^{-2t_\eps+(1-\eps)t-\frac{4(1-\eps)(t-t_\eps)^2}{4t} 
} \gamma \le \gamma 
\]
for $t\ge t_\eps$.
Since $v(t,\cdot)$ is obviously increasing on $(-\infty,-1)$, 
it follows that $v$ is a sub-solution of \eqref{1.5}  on the domain 
\begin{equation}\label{dom-D}
D\equiv ([0, t_\eps)\times\bbR) \cup \{(t,x) \,|\, t_\eps\ge t\hbox{ and } x< -2\sqrt{1-\eps}(t-t_\eps)-1  \}.
\end{equation}
Moreover,  $w$ is a solution of \eqref{1.5}, 
\[
v(0,x)= e^{-2t_\eps} \gamma \chi_{[-1,0]}(x)\le w(0,x),
\]
 and 
\[
v(t,-2\sqrt{1-\eps}(t-t_\eps)-1)\le\gamma \le w(t,-2\sqrt{1-\eps}(t-t_\eps)-1)
\]
for $t\ge t_\eps$. Thus $v\le w\le u$ on $\bar D$. Since $-2\sqrt{1-\eps}(t-t_\eps)-1\ge m-2t$ 
implies $(t,x)\in\bar D$ whenever $x\le m-2t$, the result now follows 
with $k_\eps \equiv e^{-2t_\eps}$. \hfill $\qed$

\medskip

\noindent {\bf Proof of Lemma~\ref{L.1.3}.}
Assume that
\[
u(t',x)\ge C_\eps u(0,M)\sqrt{|t'|}\, e^{\sqrt{\lambda_M  -1}\,
x +(\lambda_M   -\eps) t'}
\]
for some $t'\le -1$ and $x\in[0,c_\eps(-t'-1)]$, let $t\equiv t'+1\le 0$, and define
\[
\beta\equiv \frac{x}{2|t|\sqrt{\lambda_M-1}} \le  \frac {\lambda_M-\eps}{2(\lambda_M-1)} < 1.
\]
By the Harnack inequality and parabolic regularity that there 
exists $c_0\in(0,e^{-\lambda_M}\tht/ 2)$ (depending on $a,\tht$) such that 
\begin{equation}\label{mar231}
u(t,z) \ge c_0C_\eps u(0,M) \sqrt{|t|+1}\, e^{\sqrt{\lambda_M -1}\,x +(\lambda_M  
-\eps) t}  
\end{equation}
for all $z\in(x-1,x)$. Note that the right side of (\ref{mar231}) is below 
$\theta/2$ since $u(t,x)\le 1$.
Then Lemma~\ref{L.1.2} with $l\equiv x$ and $m\equiv 2M$ shows that for  
$y\in[0,2M]$ and $C'_\eps\equiv k_\eps c_0C_\eps$ 
(with $k_\eps$ from that lemma and using $\sqrt{\lambda_M-1}>1$) we have
\begin{align*}
u(t+\beta |t|,y) & \ge  C'_\eps u(0,M) \sqrt{|t|+1}\, 
e^{\sqrt{\lambda_M -1}\,x +(\lambda_M  -\eps) t} 
e^{(1-\eps)\beta |t|} \int_{x-1}^{x} 
\frac {e^{-{|y-z|^2} /{4\beta |t|}}}{\sqrt {4\pi\beta  |t|}}   \, dz 
\\ & \ge \frac{C'_\eps u(0,M)}{\sqrt{4\pi}} 
e^{\sqrt{\lambda_M -1}\,x +\lambda_M t  -\frac{x^2} {4\beta |t|} + \beta |t|}. 
\end{align*}
The normalization $\|\psi_M\|_\infty=1$ and the comparison principle then give
\[
u(0,z) \ge \min \left\{ \tht, e^{\lambda_M(1-\beta) |t|}  
\frac{C'_\eps u(0,M)}{\sqrt{4\pi}} e^{\sqrt{\lambda_M -1}\,x -\lambda_M |t| 
-\frac{x^2}  {4\beta |t|} + \beta |t|}  \right\}  \psi_M(z)  = 
 \min \left\{ \tht, \frac{C'_\eps  u(0,M)}{\sqrt{4\pi}} \right\} \psi_M(z)
\]
for any $z\in\bbR$.
Taking $z=M$ and $C_\eps= 4\sqrt{\pi}/k_\eps c_0\psi_M(M)$, 
it follows that 
\[
u(0,M)\ge \min\{\tht\psi_M(M),2u(0,M)\},
\]
which contradicts \eqref{1.1a} and $u(0,M)>0$. Thus, \eqref{1.2} holds for this 
$C_\eps$. 
\hfill $\qed$



\subsubsection*{The case of general inhomogeneities}
 
We now dispense with the assumptions of $a(x)-1$ compactly supported and \eqref{apr2902}. The proof of \eqref{10.2}
easily extends to the case of \eqref{10.1} and \eqref{apr3004}.
First, pick $\eps\in(0,c-2)$ (recall that $c>2$) such that $(\blambda-2\eps)/\sqrt{\blambda-1}> c $ and then $\tht>0$ such that $f(x,u)\ge (a(x)-\eps/2)u$ for $u\le\tht$. Next, 
choose $M$ large 
enough so that $a(x)\le 1+\eps$ outside $(0,2M)$ (after a shift in $x$ as 
before) and the principal eigenvalue 
$\lambda_M$ $(<\lambda-\eps/2)$ of the operator
\[
\partial_{xx}+a(x)-\eps/2
\]
on $(0,2M)$ with Dirichlet boundary conditions satisfies $\lambda_M>\lambda-\eps$. Thus $c_\eps\equiv (\lambda_M-\eps)/\sqrt{\lambda_M-1}>c$, so we  can again let $c'\equiv c_\eps>c$. 

Then Lemma \ref{L.1.2}  holds for the chosen $\eps,\tht$ without a change in the proof, even though now we have only $f(x,u)\ge (1-\eps/2)u$  
for $u\le\tht$.  Lemmas \ref{L.1.3} and \ref{lem1apr30} are also unchanged. The only change in the proof of non-existence of fronts in
Theorem \ref{T.1.1} is that one has to take
\[
v_{t_0}(t,x) \equiv C_0 u(0,M) e^{x-ct_0 + (2+\eps)(t+t_0)} + C_0 u(0,M) 
e^{4M-x+ct}.
\]
Since $c>2+\eps$, we again obtain
\[
u(t,x)\le C_2 u(0,M) e^{-|x|+ct} 
\]
for $t\le 0$ and $x\in\bbR$, so \eqref{10.2} as well as non-existence of fronts follow. 

A bump-like solution is now obtained as a limit  of solutions $u_n(t,x)$ defined on $(-n,\infty)\times\bbR$ with initial data $u(-n,x)=C_n \psi(x)$.  Here $0<C_n\to 0$ are chosen so that $u_n(0,0)=1/2$, and parabolic regularity ensures that a global in time solution $u$ of \eqref{1.1} can be obtained as a locally uniform limit on $\bbR^2$ of $u_n$, at least along a subsequence.  Since $C_{n}e^{\blambda (t-n)}\psi(x)$ is a supersolution of \eqref{1.1}, we have  $C_{n}e^{\blambda}\ge C_{n-1}$.  Since $C_{n}e^{(\blambda-\eps_n) (t-n)}\psi(x)$ is a subsolution of \eqref{1.1} on $[-n,-n+1]$ provided 
\[
\eps_n\equiv \sup_{(x,u)\in\bbR\times(0,C_ne^\blambda)} \left[a(x)-\frac{f(x,u)}u \right] \quad(\le CC_n^\delta e^{\blambda\delta} \text{ by \eqref{10.1}})
\]
and using $\|\psi\|_\infty=1$, we have  $C_{n} e^{\blambda-\eps_n}\le C_{n-1}$. Thus $C_n$ decays exponentially and then so does $\eps_n$. As a result, $C_ne^{\blambda n}\to C_\infty\in (0,\infty)$ and so $u_n(t,x)\le 2C_\infty e^{\blambda t}\psi(x)$ for all large $n$ and all $(t,x)$. Thus the limiting solution $u$ also satisfies this bound and it is therefore bump-like.

The proof of uniqueness of global solutions also extends to \eqref{apr3004}, but this time \eqref{apr2902} is necessary in order to obtain \eqref{may301} and to then apply Proposition 2.5 from~\cite{HP}.

\section{Fronts with speeds $c\in (2,\blambda/\sqrt{\blambda-1})$: The proof of Theorem~\ref{T.1.2}}

 First note that the proof of existence of bump-like solutions from Theorem \ref{T.1.1} works for any $a_->0$ and extends to $\blambda<2$, so we are left with proving existence of fronts.  
 
 Assume that $a(x)=1$ outside $[-M,M]$ and also (for now) that \eqref{apr2902} holds.
Consider any $c\in(2,\blambda/\sqrt{\blambda-1})$.  
We will  construct a positive solution $v$ and a  
sub-solution $w$ to the PDE
\[
u_t=u_{xx}+a(x)u,
\]
such that $w\le \min\{v,\tht\}$ and both move to the right with speed $c$ (in a sense to be specified later). It follows 
that $v$ and $w$ are a  supersolution and a subsolution to  (\ref{1.1}), and we will see later that this ensures the existence of a transition front $u\in(w,v)$ for \eqref{1.1}. 

For any $\gamma\in(\blambda,2)$ let $\phi_\gamma$ be the unique solution of
\beq \lb{8.1}
\phi_\gamma''+a(x)\phi_\gamma=\gamma\phi_\gamma,
\eeq
with $\phi_\gamma(x)=e^{-\sqrt{\gamma-1}\,x}$ for $x\ge M$.  We claim that then 
\beq\lb{8.2}
\phi_\gamma>0.
\eeq
Indeed, assume $\phi_\gamma(x_0)=0$ and let $\psi_\gamma$ be the solution of \eqref{8.1} with $\psi_\gamma(x)=e^{\sqrt{\gamma-1}\,x}$ for $x\ge M$. Then $\phi_\gamma-\eps\psi_\gamma$ would have at least two zeros for all small $\eps$ (near $x_0$ and at some $x_1\gg M$). Since $\gamma>\blambda =  \sup\sigma(\partial_{xx}^2+a(x))$, this would contradict the Sturm oscillation theory, so \eqref{8.2} holds.  Since there are $\alpha_\gamma,\beta_\gamma$ such that
\[
\phi_\gamma(x)=\alpha_\gamma e^{-\sqrt{\gamma-1}\,x}+\beta_\gamma e^{\sqrt{\gamma-1}\,x}
\]
for $x\le -M$, it follows  that $\alpha_\gamma>0$. 

This means that the function
\[
v(t,x) \equiv e^{\gamma t} \phi_\gamma(x) >0
\]
is a supersolution of \eqref{1.1} (if we define $f(x,u)\equiv 0$ for $u>1$).  Notice that in the domain $x>M$, the graph of $v$ moves to the right at {\it exact} speed  $\gamma/\sqrt{\gamma-1}$ as time increases.  This is essentially true also for $x\ll -M$ (since $\phi_\gamma(x)\approx \alpha_\gamma e^{-\sqrt{\gamma-1}\,x}$ there), so $v$ is a supersolution moving to the right at speed $\gamma/\sqrt{\gamma-1}$  in the sense of Remark 1 after Theorem \ref{T.1.2}.  

Next let $0<\eps'\le \eps$ and $A>0$ be large, and define
\[
w(t,x)\equiv e^{\gamma t} \phi_\gamma(x) - A e^{(\gamma+\eps) t} \phi_{\gamma+\eps'}(x).
\]
Then $w$ satisfies
\beq \lb{8.2a}
w_t=w_{xx} +a(x)w - (\eps-\eps')A e^{(\gamma+\eps) t} \phi_{\gamma+\eps'}(x).
\eeq
If we define $f(x,u)\equiv 0$ for $u<0$, then $w$ will be a subsolution of \eqref{1.1} if  $\sup_{(t,x)} w(t,x)\le \tht$, due to \eqref{apr2902}. We will now show that we can choose $\eps,\eps',A$ so that this is the case.

For large $t$ such that $\supp \, w_+ \subseteq (M,\infty)$ (namely, $t>\eps^{-1}(\sqrt{\gamma+\eps'-1}\,M-\sqrt{\gamma-1}\,M-\log A)$), the maximum  $\max_x w(t,x)$ is attained at $x$ such that
\beq\lb{8.3}
 \sqrt{\gamma-1}\, e^{\gamma t} e^{-\sqrt{\gamma-1}\,x} =  A \sqrt{\gamma+\eps'-1}\, e^{(\gamma+\eps) t} e^{-\sqrt{\gamma+\eps'-1}\,x},
\eeq
that is, at
\beq\lb{8.4}
x_t \equiv \frac {1}{\sqrt{\gamma+\eps'-1} - \sqrt{\gamma-1}} \left[ \eps t +\log \left( A \frac{\sqrt{\gamma+\eps'-1}}{\sqrt{\gamma-1}}\right) \right].
\eeq
If we define
\[
\kappa =\kappa(\eps',\gamma) \equiv \frac{\sqrt{\gamma-1}} {\sqrt{\gamma+\eps'-1} - \sqrt{\gamma-1}} >0,
\]
then we have 
\beq\lb{8.5}
w(t,x_t) = e^{(\gamma - \eps \kappa) t} A^{-\kappa} \left( \frac {\sqrt{\gamma+\eps'-1}}{\sqrt{\gamma-1}}\right)^{-\kappa-1} \left( \frac {\sqrt{\gamma+\eps'-1}}{\sqrt{\gamma-1}} -1  \right)
\eeq
for  $t\gg 1$.
So if $\eps\ge \eps'$ are chosen so that $\eps\kappa=\gamma$ (this is possible because $\gamma>2(\gamma-1)$), then  $\max_x w(t,x)$ is constant for $t\gg 1$. 

The same argument works for $t\ll -1$, with $A \alpha_{\gamma+\eps'}/\alpha_\gamma$ in place of $A$ in \eqref{8.3}---\eqref{8.5}, as well as with all three equalities holding only approximately due to the term $\beta_\gamma e^{\sqrt{\gamma-1}\,x}$.  Nevertheless, the equalities hold in the limit $t\to-\infty$, and $\max_x w(t,x)$ has a positive limit as $t\to-\infty$.  Therefore $\max_x w(t,x)$ is uniformly bounded in $t$, and this bound converges to 0 as $A\to\infty$, due to \eqref{8.5}. We can therefore pick $A$ large enough so that $\sup_{(t,x)} w(t,x)\le \tht$, so that $w$ is now a subsolution of \eqref{1.1}.
Note that $\eps\kappa=\gamma$ also implies that $x_t$ (and hence $w$) moves to the right with speed
\[
\frac {\eps}{\sqrt{\gamma+\eps'-1} - \sqrt{\gamma-1}} = \frac \gamma {\sqrt{\gamma-1}}
\] 
(in the sense of $\sup_t |x_t- {\gamma t} /{\sqrt{\gamma-1}}|<\infty$).

So given $c\in(2,\blambda/\sqrt{\blambda-1})$ let us pick $\gamma\in(\blambda,2)$  such that $c=\gamma /\sqrt{\gamma-1}$ (and then choose $\eps,\eps',A$ as above).  Then we have a subsolution $w$ and a supersolution $v$ of \eqref{1.1} with $v>\max\{w,0\}$, $\max_x w(t,x)$ bounded below and above by positive constants,  with the same decay as $x\to\infty$, and with $v\to\infty$ and $w\to -\infty$ as $x\to -\infty$.  Moreover, $v$ and $w$ are moving at the same speed $c$ to the right, in the sense that points where $\max_x w(t,x)$ is achieved and where, say, $v(t,x)=1/2$, both move to the right with speed $c$ (exact for $t\gg 1$ and almost exact for $t\ll -1$).

A standard limiting argument (see, for instance, \cite{fm}) now recovers a 
global in time solution to (\ref{1.1}) that is sandwiched between $v$ and $w$.  Indeed, we obtain it as a locally uniform limit (along a subsequence if needed) of solutions $u_n$ of \eqref{1.1} defined on $(-n,\infty)\times\bbR$, with initial condition $u_n(-n,x)\equiv \min\{ v(-n,x),1\}$, so that $u\in(\max\{w,0\},\min\{v,1\})$ by the strong maximum principle.
Another standard argument based on the global stability of the constant solution $1$ (on the set of solutions $u\in(0,1)$), same speed $c$ of $v$ and $w$, and uniform boundedness below of $\max_x w(t,x)$ in $t$ shows that $u$ has to be a transition front moving with  speed $c$, in the sense of Remark 1 after Theorem \ref{T.1.2}. 

This proves the existence-of-front part of Theorem \ref{T.1.2} when \eqref{apr2902} holds.  In that case we could even have chosen $\eps'=\eps$ so that $\eps\kappa=\gamma$ because then $\lim_{\eps\to 0} \eps\kappa = 2\sqrt{\gamma-1}<\gamma <\infty = \lim_{\eps\to \infty} \eps\kappa$. If we only have \eqref{10.1}, we need to pick $\eps'<\eps$ such that $\eps\kappa=\gamma$ and the last term in \eqref{8.2a} to be larger than $Cw(t,x)^{1+\delta}$ where $w(t,x)>0$, so that $w$ stays a subsolution of \eqref{1.1}.  For the latter it is sufficient if
\beq \lb{8.7}
(\eps-\eps')A e^{(\gamma+\eps) t} e^{-\sqrt{\gamma+\eps'-1}\,x}\ge C_1 e^{-(1+\delta)\sqrt{\gamma-1}\,x}
\eeq
where $w(t,x)>0$, with some large $C_1$ depending on $C$, $\phi_\gamma$, $\phi_{\gamma+\eps'}$.  If we let $y\equiv x-ct = x-\gamma t/\sqrt{\gamma -1}$ and use $\eps\kappa=\gamma$, this boils down to 
\beq \lb{8.8}
\sqrt{\gamma+\eps'-1}\,y<(1+\delta)\sqrt{\gamma-1}\,y + \log \frac{(\eps-\eps')A}{C_1}
\eeq
when $w(t,ct+y)>0$.  Notice that for say $A=1$, the leftmost point where $w(x,t)=0$ stays uniformly (in $t$) close  to $ct$ (say distance $d(t)\le d_0$), and only moves to the right if we increase $A$.  Therefore we only need to pick $\eps'<\eps$ such that $\sqrt{\gamma+\eps'-1}\le (1+\delta)\sqrt{\gamma-1}$ and $\eps\kappa=\gamma$, and then $A>1$ large enough so that \eqref{8.8} holds for any $y\ge -d_0$.  The rest of the proof is unchanged.
\hfill $\qed$

\section{Nonexistence of fronts with speeds $c>\blambda/{\sqrt{\blambda-1}}$: 
The proof of Theorem \ref{T.1.3}}

Assume $a(x)\equiv 1$ outside $[-M_0,M_0]$ and let us denote the roots of $ r^2-cr+1=0$ by
\[
r_\pm(c)=\frac{c\pm\sqrt{c^2-4}}{2}.
\]
Notice that if $\blambda\le 2$ and $c>\blambda/{\sqrt{\blambda-1}}$, then 
\beq \lb{9.1}
0< r_-(c) < \sqrt{\blambda-1} \qquad \text{and}\qquad r_+(c)> \frac1{\sqrt{\blambda-1}}.
\eeq
Also recall that we denote by $X(t)$ the  right-most point $x$ such that $u(t,x)=1/2$.
The proof of Theorem \ref{T.1.3} relies on the following upper and lower exponential bounds on the solution ahead of the front (at $x\ge X(t)$). 

\begin{lemma}
\label{l3.1}
Let $c>2$ and $u(t,x)$ be a transition front for \eqref{1.1} moving with 
speed $c$. Then for any $\e>0$ there exists $C_\eps>0$ such that
\begin{equation}\label{apr2201}
u(t,x)\leq C_\eps e^{-(r_-(c)-\e)(x-X(t))}\ \ \ \hbox{for $x\geq X(t)$.}
\end{equation}
\end{lemma}

\begin{lemma}
\label{l3.2}
Assume that the function $a(x)$ is even and that \eqref{apr2904} has a unique eigenvalue $\blambda > 1$.
Let $c> \blambda/\sqrt{\blambda-1}$ and $u(t,x)$ be a transition front for \eqref{1.1} moving with 
speed $c$. 
Then for all $\e>0$, there is $C_\e>0$ and $T > 0$ such that:
$$
u(t,x)\geq C_\e e^{-(r_-(c)+\e)(x-X(t))}\ \ \ 
\hbox{for $t \geq T$ and $x\geq X(t)$.}
$$
\end{lemma}

\noindent{\bf Proof of Theorem \ref{T.1.3}.}
Let us assume $\blambda\in(1,2]$ since the cas $\blambda>2$ has already been proved in Theorem \ref{T.1.1}.
 Assume that there exists a  transition  front $u(t,x)$ with speed 
\begin{equation}\label{apr1906}
c>\blambda/{\sqrt{\blambda-1}}.
\end{equation} 
We first wish to prove the following estimate: for all 
$\e>0$,
there exists $C_\e>0$ such that
\begin{equation}
\label{e6.1}
u(t,x)\leq C_\e e^{(\blambda-\e)t-\sqrt{\blambda-\e-1}\,x}~~\hbox{for all $x\geq 0$ and $t\le 0$.}
\end{equation}
From Lemma~\ref{L.1.3}, the estimate is true for $x=0$ and, more generally,
on every bounded subset of $\RR_+$, so let us extend it to the
whole half-line. For this, we notice that, for all $t\leq 0$, we have
\begin{equation}
\label{e6.2}
u(t,x)\leq Ce^{t},~~\hbox{for $x\ge 0$}.
\end{equation}
Indeed, the function
$$
\alpha(t)=\int_{M_0}^{+\infty}u(t,x)\ dx,
$$
which is finite due to Lemma \ref{l3.1},
solves
$$
\alpha'-\alpha=-u_x(t,M_0)-\int_{M_0}^{+\infty}(u(t,x) - f(u(t,x))\ dx.
$$
From parabolic regularity and \eqref{e6.1} for $x$ on compact intervals, we have
$|u_x(t,M_0)|\le Ce^{(\blambda-\e)t}$ for $t\le 0$. From Lemma \ref{l3.1}, the fact that $u$ travels with a positive speed, and $a(x)=1$ for $x\ge M_0$, we 
have $f(u(t,x))=u(t,x)$ for $x\ge M_0$ and $t\ll -1$.
Hence we have 
$$
\alpha'-
\alpha=O(e^{(\blambda-\e)t})
$$
for $t\ll -1$,
which implies $\alpha(t)=O(e^{t})$ for $t\leq0$ since $\blambda>1$. Estimate 
\eqref{e6.2} then follows from parabolic regularity.

Then, we set
$$
w(t,x)=e^{-t}u(t,x)-C_\e e^{(\blambda-\e-1)t-\sqrt{\blambda-\e-1}(x-M-1)}.
$$
Since (\ref{e6.1}) holds on compact subsets of $\RR_+$, we have
\begin{eqnarray*}
&&w_t-w_{xx}\leq 0\hbox{ for $t\leq0$, $x\geq M_0$},\\
&&w(t,M_0)\leq 0\hbox{ for $t\leq0$.}
\end{eqnarray*}
From \eqref{e6.2} (and $\blambda>1$) the function $w$ is
bounded on $\RR_-\times[M_0,+\infty)$. Consequently, it cannot attain
a positive maximum, and there cannot be a sequence $(t_n,x_n)$ such that
$w(t_n,x_n)$ tends to a positive supremum. This implies that
$w$ is negative, hence estimate \eqref{e6.1} for $x\ge M_0$ follows. It also
holds on $[0,M_0]$ due to parabolic regularity. 

Let us now turn to positive times. 
The function $v(t,x)=u(t,x+ct)$ solves  
\begin{eqnarray*}
&&v_t-v_{xx}-cv_x\leq v\hbox{ for $t\geq0,\ x\geq M_0$,}\\
&&v(t,M_0)\leq 1\hbox{ for $t\geq0$},\\
&&v(0,x)\leq C_\eps e^{-\sqrt{\blambda-1-\e}\,x},
\end{eqnarray*}
the last inequality due to \eqref{e6.1}.
Since for small enough $\e>0$ we have $r_-(c)<\sqrt{\blambda-\e-1}<r_+(c)$, the stationary function $e^{-\sqrt{\blambda-1-\e}\,x}$ is a super-solution to 
$$
v_t-v_{xx}-cv _x=v.
$$
This in turn implies
$
v(t,x)\leq C_\eps e^{-\sqrt{\blambda-1-\e}\,x}$ for small $\eps>0$. Using the fact that the front travels with speed $c$, we get
$$
u(t,x)\leq C e^{-\sqrt{\blambda-1-2\e}(x-X(t))}
$$
with a new $C$. This contradicts Lemma \ref{l3.2} since $r_-(c)<\sqrt{\blambda-1}$.
\hfill $\qed$
\medskip

The rest of the paper contains the proofs of Lemmas~\ref{l3.1} and~\ref{l3.2}.

\section{An upper bound for fronts with speed $c > \blambda / \sqrt{\blambda - 1}$: The proof of Lemma~\ref{l3.1}}
 
 It is obviously sufficient to prove that
for any $\eps>0$ there exists $x_\eps$  such that for any $t\in\RR$ we have 
\begin{equation}\label{apr1603_0}
u(t,x)\leq   e^{-(r_-(c)-\e)(x-X(t))}\ \ \ \hbox{for $x\geq X(t)+x_\eps$}.
\end{equation}
Therefore assume, towards contradiction, that there exists
$\e>0$ and $T_n\in\bbR$, $x_n\to +\infty$ such that 
$$
u(T_n,X(T_n)+x_n)\geq e^{-(r_-(c)-\e)x_n}.
$$
By the Harnack 
inequality, there is a constant $\delta>0$ such that 
\begin{equation}\label{apr1604}
u(T_n-1,X(T_n)+x)\geq \delta e^{-(r_-(c)-\e)x_n}\ \ \ \hbox{for $x\in[x_n,x_n+1]$}.
\end{equation}

As $u$ satisfies (\ref{apr2816}) and moves with speed $c$, we know that for every 
$\galpha>0$ we have
\[
\lim_{s\to+\infty}\sup_{T\in\bbR, \,\,x\ge X(T)+(c+\galpha)s}u(T+s,x)=0.
\]
Therefore, for every  $\galpha  >0$ there is $x_\alpha>0$ 
such that for any $T\in\bbR$,
$$
f(u(t,x)) =  u(t,x) \ \ \ \hbox{for  $t\geq T$ and $x\geq X(T)+(c+\galpha  )(t-T)+x_\alpha$}
$$ 
Then from $u\le 1$ we have for $t\ge T$
$$
u_t-u_{xx} =   a(x) u + a(x)\left(f(u) - u\right) \geq u-C\un_{x\leq X(T)+(c+\galpha  )(t-T)+x_\alpha}
$$
with $C=\|a\|_\infty$. Thus we have 
$$
u(t,x)\geq\di{e^{t}\int_{\Rm}\frac{e^{-\frac{(x-y)^2}{4(t-T)}}}{\sqrt{4\pi 
(t-T)}}u(T,y)\ dy-
C\int_T^t\int_{-\infty}^{x_\alpha+(c+\galpha  )s}\frac{e^{-\frac{(x-y)^2}{4(t-s)}+(t-s)} 
}{\sqrt{4\pi (t-s)}}\ dyds}
=:I(t,x)-I\!I(t,x)
$$
We are going to evaluate $I(t,x)$ and $I\!I(t,x)$ for $T=T_n-1$ at 
$$
(t,x)=(t_n,z_n):= \left(T_n-1+\frac{x_n}{\sqrt{c^2-4}}, X(T_n)+\frac{cx_n}{\sqrt{c^2-4}} \right),
$$
and show that $I(t_n,z_n)\to+\infty$ faster than $I\!I(t_n,z_n)$ provided $\alpha>0$ is small enough, giving a 
contradiction with $u(t,x)\le 1$.  

Fix $n$ and for the sake of simplicity assume $T_n=1$ and $X(T_n)=0$ (this can be achieved by a translation in space and time). So $T=0$ and by  (\ref{apr1604}) we have
\begin{eqnarray*}
I(t_n,z_n)\geq  e^{t_n}\int_{x_n}^{x_n+1}\frac{e^{-\frac{(z_n-y)^2}{4t_n}}}{\sqrt{4\pi 
t_n}}u(0,y)\ dy
\geq\frac{\delta}{\sqrt{4\pi t_n}}e^{t_n-(r_-(c)-\e)x_n}\int_0^1e^{-\frac{(z_n-x_n-z)^2}{4t_n}}\ 
dz.
\end{eqnarray*}
Note that for $z\in[0,1]$ we have
\[
\di{\frac{(z_n-x_n-z)^2}{t_n}=\frac{(z_n-x_n)^2}{t_n}}+O(1),
\] 
thus with some $n$-independent $q>0$ we have
$$
I(t_n,z_n)\geq\frac{q\delta}{\sqrt{4\pi t_n}}e^{-\frac{(z_n-x_n)^2}{4t_n}+t_n-(r_-(c)-\e)x_n}.
$$
The exponent is easily evaluated using the relations 
$x_n=\sqrt{c^2-4}\, t_n$, $z_n-x_n=2r_-(c)t_n$, and $r_-(c)^2+ \sqrt{c^2-4}\,r_-(c) -1=0$, leading to
\beq
\label{e3.1}
I(t_n,z_n)\geq\frac{q\delta}{\sqrt{4\pi t_n}}e^{(\e\sqrt{c^2-4}-\alpha)t_n}.
\eeq

To estimate $I\!I(t_n,z_n)$, notice that we have (using $z_n=ct_n$ and with $z:=y-z_n$)
\begin{eqnarray*}
&&I\!I(t_n,z_n)\le C{\int_0^{t_n}\int_{-\infty}^{x_\alpha+(c+\galpha  )s-z_n}
\frac{e^{t_n-s-\frac{z^2}{4(t_n-s)}}}{\sqrt{4\pi(t_n-s)}}dzds}\\
&&~~~~~~~~~~~~
=C\di{\int_0^{t_n}\biggl(\int_{x_\alpha-c(t_n-s)}^{ x_\alpha-c(t_n-s)+\galpha   s}+\int_{-\infty}^{x_\alpha-c(t_n-s)}\biggl)\ \frac{e^{t_n-s-\frac{z^2}{4(t_n-s)}}}{\sqrt{4\pi(t_n-s)}}
dzds } \\
&&~~~~~~~~~~~~
=: I\!I_1(t_n,z_n)+I\!I_2(t_n,z_n)  .         
\end{eqnarray*}
Using the estimate
$$
\int_{-\infty}^{x_\alpha-c(t_n-s)} \frac{e^{-\frac{z^2}{4(t_n-s)}}}{\sqrt{t_n-s}}  {dz}\leq 
C_\alpha\frac{e^{-\frac{c^2(t_n-s)}4}}{\sqrt{t_n-s}}
$$
and $c>2$, we have $I\!I_2(t_n,z_n)=O(1)$ as ${n\to+\infty}$. In order to estimate $I\!I_1(t_n,z_n)$, we represent $\zeta:=z+c(t_n-s)\in [x_\alpha, x_\alpha+\galpha  s]$ so that
\[
t_n-s-\frac{z^2}{4(t_n-s)}=t_n-s-\frac{c^2(t_n-s)^2+\zeta^2-2c(t_n-s)\zeta}{4(t_n-s)}\le
\farc{c\zeta}{2}\le cx_\alpha+c\galpha  t_n.
\]
It follows that
$$
I\!I_{1}(t_n,z_n)\leq \galpha  t_ne^{cx_\alpha+c\galpha  t_n}\int_0^{t_n} 
\frac{ds}{\sqrt{4\pi(t_n-s)}}\leq C\galpha  t_n^{3/2}
e^{{cx_\alpha+c\galpha  }t_n}\le C_\alpha  e^{{2c\galpha  }t_n}.
$$
We now choose $\alpha>0$ so that $\eps\sqrt{c^2-4}-\alpha>2 c\galpha  $.  
Using \eqref{e3.1}, it follows that $u(t_n,z_n)=I(t_n,z_n)-I\!I(t_n,z_n)>1$ for all large $n$, a contradiction.
This finishes the proof of Lemma~\ref{l3.1}.  \hfill $\qed$

\section{A lower bound for fronts with speed $c > \blambda / \sqrt{\blambda - 1}$: The proof of Lemma \ref{l3.2}}
\subsection{A heat kernel estimate}
We will need a rather precise information on the behavior, for 
large $x$ and $t$, of the solutions of the Cauchy problem
\begin{eqnarray}
\label{e4.1}
&&u_t-u_{xx}-A(x)u=0,~~t>0,\ x\in\RR,\\
&&u(0,x)=u_0(x).\nonumber
\end{eqnarray}
The function $B(x)=A(x)-1$ is assumed to be nonnegative and to have compact 
support, in an interval $[L-M_0,L+M_0]$. 
Basically, $A$ should be thought of as a translate
of the function $a$: in the proof of Lemma~\ref{l3.2} below, the number $M_0$ will be of fixed size, 
the number $L$ will vary arbitrarily.
A lot -- most probably, including our estimate below -- is known about solutions
of (\ref{e4.1}). See, for instance, \cite{Pinch} and the references therein. 
However we were not able
to find in the literature an estimate of the type \eqref{e4.2} below. Moreover, 
the proof is short, so it is
worth presenting it in reasonable detail.
Denote by $G(t,x,y)$ the heat kernel of \eqref{e4.1}, i.e. the function such 
that the solution $u(t,x)$ is
$$
u(t,x)=\int_{-\infty}^{+\infty}G(t,x,y)u_0(y)\ dy.
$$
Let us also denote by $H(t,z)$ the standard heat kernel:
$$
H(t,z)=\frac{e^{-z^2/4t}}{\sqrt{4\pi t}}.
$$
\begin{proposition}
\label{p4.1}
Assume the function $B(x-L)$ to be even and nonnegative, and that the eigenvalue problem
$$
\phi_0''+(1+B(x-L))\phi_0 =\lambda \phi_0
$$
has a unique eigenvalue $\lambda > 1$. Let $\phi_0>0$ be the eigenfunction with $\|\phi_0\|_2=1$. Then 
we have
\begin{equation}
\label{e4.3}
G(t,x,y)\geq e^tH(t,x-y).
\end{equation}
for all $x,y \in \RR$. Conversely, if $x<L-M_0$ and $y>L+M_0$, or $y<L-M_0$ and $x>L+M_0$, then
there is a smooth function $\psi_0$ such that $\psi_0(x)=O(e^{-\sqrt{\lambda-1}\vert x\vert})$ for $ \vert x-L\vert\geq 2M_0$, and such that, for all $\e>0$ we have
\begin{equation}
\label{e4.2}
\vert G(t,x,y) - \biggl(e^{\lambda t}\phi_0(x)\phi_0(y)+e^t(H(t,.)*\psi_0)(x-y)\biggl)\vert\leq Ce^{t+C\vert x-y\vert/t}H(t,x-y).
\end{equation}
Also, there is $C>0$, depending on $M_0$ but not on $L$, such that if $x,y<L-M_0$ or $x,y>L+M_0$, we 
have
\begin{equation}
\label{e4.4}
G(t,x,y) - \biggl(e^{\lambda t}\phi_0(x)\phi_0(y)+e^t(H(t,.)*\psi_0)(x+y-2L)\biggl) \leq 
Ce^{t+C\vert x+y-2L\vert/t}H(t,x+y-2L).
\end{equation}
\end{proposition}
\noindent{\bf Proof.} The lower bound \eqref{e4.3} is obvious, because $A(x)\geq 1$. So, let us examine
the upper bound. First, we may without loss of generality assume $L=0$, the 
result will just follow by translating $x$ and $y$ by the amount $L$. Also, it is enough to replace $A(x)$ 
by $B(x)$ (thus we deal with a compactly supported potential), at the
expense of multiplying the final result by $e^t$. Our proof will use some basic 
facts of eigenfunction expansions, see \cite{Lev},
that we recall now. For $k\in\RR^*$, let us denote by $f(x,k)$ the solution of 
\begin{equation}
\label{e4.5}
-\phi''=(B(x)+k^2)\phi,\ \ \ x\in\RR
\end{equation}
satisfying
\begin{equation}
\label{e4.6}
f(x,k)=e^{ikx}\ \ \ \hbox{for $x\geq M_0$}
\end{equation}
and let us denote by $g(x,k)$ the solution of \eqref{e4.5}
such that 
\begin{equation}
\label{e4.7}
g(x,k)=e^{-ikx}\ \ \ \hbox{for $x\leq -M_0$.}
\end{equation}
Denoting by $W(u(x),v(x))$ the Wronskian of two solutions $u$ and
$v$ of \eqref{e4.5}, let us set
$$
a(k)=-\frac1{2ik}W(f(x,k),g(x,k)),\ \ \ \ b(k)=\frac1{2ik}W(f(x,k),g(x,-k))
$$
and
\begin{equation}
\label{e4.14}
c(k)=-b(-k),\ \ \ d(k)=a(k).
\end{equation}
We have
\begin{equation}
\label{e4.8}
\begin{array}{rll}
f(x,k)=&a(k)g(x,-k)+b(k)g(x,k)\\
g(x,k)=&c(k)f(x,k)+d(k)f(x,-k),
\end{array}
\end{equation}
and $|a(k)|^2 = 1 + |b(k)|^2$, $b(-k) = \overline{b(k)}$, and $a(-k) = \overline{a(k)}$. The following decompositions hold:
\begin{equation}
\label{e4.10}
\delta({x-y})= \phi_0(x)\phi_0(y)+\frac1{2\pi}\int_{-\infty}^{+\infty}f(x,k)\overline{f(y,k)}\ 
dk
-\frac1{2\pi}\int_{-\infty}^{+\infty}f(x,k) {f(y,k)}\frac{b(-k)}{a(k)}\ dk,
\end{equation}
and
\begin{equation}
\label{e4.11}
\delta({x-y})=\phi_0(x)\phi_0(y)+\frac1{2\pi}\int_{-\infty}^{+\infty}g(x,k)\overline{g(y,k)}\ 
dk
+\frac1{2\pi}\int_{-\infty}^{+\infty}g(x,k) {g(y,k)}\frac{b(k)}{a(k)}\ dk.
\end{equation}
These decompositions may also be viewed as a consequence of Agmon's limiting 
absorption
principle, see \cite{Ag}, Theorem 4.1. Consequently, we have the representation
\begin{equation}
\label{e4.17}
\begin{array}{rll}
G(t,x,y)=&\di{e^{(\lambda-1)t}\phi_0(x)\phi_0(y)+\frac1{2\pi}\int_{-\infty}^{+\infty}e^{-tk^2}f(x,k)\overline{f(y,k)}\ 
dk}\\
&\di{-\frac1{2\pi}\int_{-\infty}^{+\infty}e^{-tk^2}f(x,k) {f(y,k)}\frac{b(-k)}{a(k)}\ 
dk}\\
=&e^{(\lambda-1) t}\phi_0(x)\phi_0(y)+\di{\frac1{2\pi}\int_{-\infty}^{+\infty}e^{-tk^2}g(x,k)\overline{g(y,k)}\ 
dk}\\
&\di{+\frac1{2\pi}\int_{-\infty}^{+\infty}e^{-tk^2}g(x,k) {g(y,k)}\frac{b(k)}{a(k)}\ 
dk}.
\end{array}
\end{equation}
Now we prove \eqref{e4.2}. If $y<-M_0$ and $x>M_0$, the identity (\ref{e4.8}) and the first equality in (\ref{e4.17}) implies that
\begin{equation}
\label{Grepa}
\begin{array}{rll}
G(t,x,y)=&\di{e^{(\lambda-1)t}\phi_0(x)\phi_0(y)+\frac1{2\pi}\int_{-\infty}^{+\infty} \frac{e^{-tk^2}}{a(-k)}e^{ik(x-y)}\ dk}\\
=&\di{e^{(\lambda-1)t}\phi_0(x)\phi_0(y)+ (H(t,\cdot) *F_1)(x-y)}
\end{array}
\end{equation}
where $F_1$ is the inverse Fourier Transform of $\frac{1}{a(-k)}$. By using the second equality in (\ref{e4.17}), we see that the same holds for $y > M_0$ and $x < - M_0$. This function $F_1$ may be estimated by \eqref{e4.10} and \eqref{e4.8} if $y<-M_0$ and $x>M_0$:
\begin{eqnarray}
&& - \phi_0(x)\phi_0(y) = \frac1{2\pi} \int_{-\infty}^{+\infty} f(x,k) \left( \overline{f(y,k)} - f(y,k)\frac{b(-k)}{a(k)} \right) dk \no \\ 
&& \quad \quad \quad = \frac1{2\pi} \int_{-\infty}^{+\infty} \left(a(k) e^{ikx} + b(k) e^{-ikx} \right) \left( e^{-iky} - e^{iky} \frac{b(-k)}{a(k)} \right) dk \no\\
&& \quad \quad \quad = \frac1{2\pi} \int_{-\infty}^{+\infty} \frac{|a(k)|^2 - |b(k)|^2}{a(-k)} e^{ik(x- y)} dk \no \\
&& \quad \quad \quad = \frac1{2\pi} \int_{-\infty}^{+\infty} \frac{e^{ik(x- y)} }{a(-k)} dk. \no
\end{eqnarray}
The same is true for $y<-M_0$ and $x>M_0$, one just has to use \eqref{e4.11} and \eqref{e4.8}. Therefore,
\begin{equation}
\label{e7.11}
F_1=\psi_0+T_0,
\end{equation}
where $\psi_0(x)=c_0e^{-\sqrt{\lambda-1}\vert x\vert}$ for $ \vert x\vert\geq 2M_0$, $T_0$ is a compactly supported distribution, and where we have made the abuse of notation consisting in using the argument $x$ in a distribution. Combining this with (\ref{Grepa}) we obtain
$$
\begin{array}{rll}
G(t,x,y)=& e^{(\lambda-1)t}\phi_0(x)\phi_0(y) + (H(t,.)*\psi_0)(x-y)+(H(t,.)*T_0)(x-y)
\end{array}
$$
and estimate \eqref{e4.2} is concluded by a standard distributional computation. 
Now we prove \eqref{e4.4}. If $x$ and $y$ are on the same side, say $x \geq M$ and $y \geq M$, then (\ref{e4.17}) implies
\begin{equation}
\label{Grepsameside}
\begin{array}{rll}
G(t,x,y)=&e^{(\lambda-1)t}\phi_0(x)\phi_0(y) + \di\frac1{2\pi}\int_{-\infty}^{+\infty}e^{-tk^2+ik(x-y)} \ dk-\di\frac1{2\pi}\int_{-\infty}^{+\infty}e^{-tk^2+ik(x+y)}\frac{b(-k)}{a(k)} \ dk\\
=&e^{(\lambda-1)t}\phi_0(x)\phi_0(y) + H(t,x-y) +(H(t, \cdot)* F_2)(x + y), \quad \text{for}\;\;x\geq M, \;\;y\geq M
\end{array}
\end{equation}
where $F_2$ is the inverse Fourier transform of the function $ {b(-k)}/{a(k)}$. Similarly, 
\begin{equation}
\label{Grepsameside2}
\begin{array}{rll}
G(t,x,y)=&e^{(\lambda-1)t}\phi_0(x)\phi_0(y) + \di\frac1{2\pi}\int_{-\infty}^{+\infty}e^{-tk^2-ik(x-y)} \ dk-\di\frac1{2\pi}\int_{-\infty}^{+\infty}e^{-tk^2-ik(x+y)}\frac{b(k)}{a(k)} \ dk\\
=&e^{(\lambda-1)t}\phi_0(x)\phi_0(y) + H(t,x-y) +(H(t, \cdot)* F_3)(x + y), \quad \text{for}\;\;x\leq -M, \;\;y\leq -M,
\end{array}
\end{equation}
where $F_3$ is the Fourier transform of the function $ {b(k)}/{a(k)}$. It follows from \cite{Lev}, that $F_2$ and $F_3$
are $W^{1,1}$ functions. From the relations \eqref{e4.8} and decomposition \eqref{e4.10}, we find that
\begin{equation}
\label{e7.12}
F_2(x + y) = \frac1{2\pi}\int_{-\infty}^{+\infty}e^{ik(x+y)}\frac{b(-k)}{a(k)} \ 
dk = \phi_0(x)\phi_0(y) , \quad \text{for}\;\; x \geq M_0, \;\; y \geq M_0.
\end{equation}
Consequently, $F_2(z) = c_1 e^{-\sqrt{\lambda - 1}\abs{z}}$ for $z > 2M_0$. In the same fashion we have, from the decomposition \eqref{e4.11},
\begin{equation}
\label{e4.16}
F_3(x + y) = -\frac{1}{2\pi}\int_{-\infty}^{+\infty}e^{-ik(x+y)} \frac{b(k)}{a(k)}\ dk = \phi_0(x) \phi_0(y) \ \ \ \ \hbox{for 
$x\leq -M_0$, $y\leq -M_0$.}
\end{equation}
From the evenness of $B$ and the relations \eqref{e4.14}, the function $b(k)$ is 
purely imaginary, so $b(-k)/a(k) = \overline{b(k)}/a(k) = - b(k)/a(k)$. Thus, $F_3(z) = F_2(-z)$. 
And so, similarly to \eqref{e7.11} there holds
$$
F_i=\psi_0+T_i,\ \ \ \ i\in{2,3}
$$
where $T_2$ and $T_3$ are $W^{1,1}$ functions supported in $(-\infty,2M_0)$ and $(-2M_0,\infty)$, respectively. So,
for $x\geq M_0$ and $y\geq M_0$, estimate \eqref{e4.4} now follows from \eqref{Grepsameside}, since
\[
\left| (H(t,\cdot)*T_2)(x+y)\right| = \left| \int_{-\infty}^{2M_0} H(t,x + y - z)T_2(z) \,dz \right| \leq H(t,x + y-2M_0) \Vert{T_2}\Vert_{1}\]
The same argument is valid for $x\leq-M_0$ 
and $y\leq-M_0$ using \eqref{Grepsameside2}.
\hfill $\qed$

 Proposition \ref{p4.1} admits the following corollary, which takes care of what 
happens when $y$
 is in the support of $B$.
 \begin{corollary}
 \label{c4.1}
Let $\psi_0$ be defined as in Proposition \ref{p4.1} .
There is a constant $C$ such that if $y\in[L-M_0,L+M_0]$ and $x\notin[L-M_0,L+M_0]$, we have
 \begin{equation}
\label{e4.100}
 G(t,x,y) - \biggl(e^{\lambda t}\phi_0(x)\phi_0(y)+(e^tH(t,.)*\psi_0)(x-L)\biggl) \leq  
Ce^{t+C\vert x-L\vert/t}H(t,x-L).
\end{equation}
\end{corollary}
The proof is similar to that of the proposition, and is omitted.
\subsection{Proof of Lemma~\ref{l3.2}}
 
Assume the conclusion of Lemma~\ref{l3.2} to be false. Then there exists a sequence $T_n \to +\infty$,
and a sequence
$x_n\to +\infty$ such that 
\begin{equation}\label{apr1801}
u(T_n,X(T_n)+x_n)\leq e^{-(r_-(c)+\e)x_n}.
\end{equation}
\subsubsection*{ Extending \eqref{apr1801} to a large interval} 

We are going to apply the Harnack 
inequality in the following way:
if $u(t,x)$ is a global solution (in time and space) of a linear parabolic 
equation on 
$(t,x)\in\bbR\times\bbR$, there exists a universal constant $\rho\in(0,1)$ such 
that
$$
 \ u(t,x)\geq \rho u(t-1,x+\xi),~~\hbox{for all $t,x\in\Rm$ and all 
$\xi\in[-1,1]$}.
$$
Thus, for all $\xi\in[-1,1]$ and all $t\in\Rm$ and $x\in\Rm$,
and any non-negative integer $p\in{\mathbb N}$ we have
\begin{equation}\label{apr1802}
u(t,x)\geq \rho^pu(t-p,x+p\xi).
\end{equation}
Then, assumption (\ref{apr1801}) on $u$ together with (\ref{apr1802}) translate 
into 
\begin{eqnarray}\label{apr1803}
&& u(T_n-p,X(T_n)+x_n+p\xi)\leq  \rho^{-p} e^{-(r_-(c)+\e)x_n} \\
&&~~~~~~~~~~~~~~~~
=  \rho^{-p}e^{-(r_-(c)+\e)[X(T_n-p)-X(T_n)]}e^{-(r_-(c)+\e)(x_n-[X(T_n-p)-X(T_n)])},
\nonumber
\end{eqnarray}
for all $\xi\in[-1,1]$.
Note that, as $u(t,x)$ is a front moving with the speed $c$,  
there exists a constant $B>0$ so that 
\begin{equation}\label{apr1804}
X(T_n)-2c(p+B)\leq X(T_n-p)\leq X(T_n)+\frac{c}{2}(-p+B).
\end{equation}

We are going to choose $p$ as a small fraction of
$x_n$, that is, $p=[\eta x_n]$ where $[x]$ denotes the integer part of 
$x$, and $\eta>0$ is small.  Then, for any 
$x\in [(1-\eta)x_n,(1+\eta)x_n]$ we rewrite (\ref{apr1803}), 
using also (\ref{apr1804}) as
\begin{eqnarray*}
&&\!\!\!\!\!\!\!
u(T_n-p,X(T_n)+x)\leq  \rho^{-p}e^{-(r_-(c)+\e)[X(T_n-p)-X(T_n)]}
e^{-(r_-(c)+\e)(x-[X(T_n-p)-X(T_n)])+(r_-(c)+\eps)(x-x_n)}\\
&&~~~~~~~~~~~
\leq  C\rho^{-p}e^{2c(r_-(c)+\e)(p+B)}e^{-(r_-(c)+\e)(x-[X(T_n-p)-X(T_n)])+
(r_-(c)+\eps)p}
\\
&&~~~~~~~~~~~\leq C \exp\left[\left(-r_-(c)-\e+\frac{Kp}{x-[X(T_n-p)-X(T_n)]}\right)
(x-[X(T_n-p)-X(T_n)])\right],
\end{eqnarray*}
with a constant $K$ that depends on $c$, $\rho$ and $B$ but not on $p$ or $x$.
As $p=[\eta x_n]$, $x_n\to+\infty$, and $X(T_n-p)\le X(T_n)+cB/2$, choosing 
$\eta=\eps/(1+2K)$  so that 
$K\eta/(1-\eta)<\eps/2$  ensures that
$$
\frac{Kp}{x-[X(T_n-p)-X(T_n)]}\leq\frac\e2 ~~~~
\hbox{for all $x\in[(1-q\eps)x_n,(1+q\eps)x_n]$} ,
$$
for $n$ large enough.  Here we have set $q=1/(1+2K)$.

Let us now shift the origin of time and space placing it at 
$(t,x)=(T_n-p,X(T_n-p))$. And thus, in the
new coordinates we have
\beq
\label{e3.3}
u_0(x):=u(0,x)\leq C e^{-(r_-(c)+ \e/2)x}\ \ \ \ 
\hbox{for}\ x\in[(1-q\eps)x_n,(1+q\eps)x_n].
\eeq
The support of $a-1$ is also shifted accordingly: it is supported
in an interval $[L-M_0,L+M_0]$, with $L=-X(T_n-p) < - M_0$ for large $n$. 

\subsubsection*{Reduction of $u(t,x)$} 

We start from
$$
u(t,x)=\calS_a(t)u_0(x)-\di\int_0^t\calS_a(t-s)a (u -f(u))\ ds \leq 
\calS_a(t)u_0(x)-\di\int_0^t\calS_1(t-s)a(u - f(u))\ ds.
$$
and we are going to evaluate it for a well chosen $(t,x)\in\Rm_+\times\Rm_+$.
Here ${\cal S}_a$ denotes the semi-groups generated by the operator
$\partial_{xx}^2+a(x)$, and ${\cal S}_1$ is the semigroup 
generated by the operator
$\partial_{xx}^2+1$, with $a(x)$ appropriately shifted to our new coordinate 
frame.
Because $x>0$, it is outside of
$\supp(a-1)=[L-M_0,L+M_0]$; we 
 will use Proposition \ref{p4.1} and Corollary \ref{c4.1} to deal with $\calS_a(t)u_0(x)$. We have
\begin{eqnarray}\label{sep0304}
&&{\cal S}_a(t)u_0(x)\le e^t\int H(t,x-y)\biggl((u_0*\psi_0)(y)+Ce^{C\vert x-y\vert/t}u_0(y)\biggl)\ dy\nonumber\\
&&~~~~~~~~~~~~~~~+e^t\int E(t,x,y))\biggl((u_0*\psi_0)(y)+Ce^{C\vert x-y\vert/t}u_0(y)\biggl)\ dy+
e^{\blambda t}\langle\phi_0,u_0\rangle\phi_0(x)\nonumber\\
&&~~~~~~~~~~~~~=u_1(t,x)+u_2(t,x)+u_3(x),
\end{eqnarray}
where
 $E(t,x,y)=0$ if $y < L - M_0$ (since $x > L + M_0$), while
\[
E(t,x,y)= C\frac{e^{-|x+y-2L|^2/(t+1)}}{\sqrt{4\pi (t+1)}}
\]
if $y > L - M_0$. We will also set
\begin{equation}\label{oct402}
u_4(t,x)=\int_0^t\calS_1(t-s)a(u - f(u))\ ds.
\end{equation}
We will estimate each of $u_1$, $u_2$, $u_3$ and $u_4$ separately at an appropriately 
chosen point $(t_n,z_n)$ and show that $u_4$ is much larger than $u_1+u_2+u_3$ giving a contradiction.

\subsubsection*{Estimate of $u_1(t,x)$}
This is the most involved, the estimates of $u_2$ and $_3$ being simpler or similar.
First, we anticipate that $u_1$ will be evaluated at a point $(t,x)$ such that $t$ and $x$
are both large, and $x$ and $t$ of the same order of magnitude. Also, in the integral
expressing $u_1$, the integrands will be maximized at points $y$ such that $\vert x-y\vert$ is
of orde $t$. Hence, from standard convolutions
between exponentials (and the fact that $r_-(c)<\sqrt{\lambda-1}$, we do not lose any generality if we
assume the existence of a function $w_0(x)$ and a constant $C>0$ such that

(i). the function $w_0$ is bounded on $\RR$,

(ii). there is a constant $C>0$ such that (even if it means restricting $q$ a little)
$$
\begin{array}{rll}
&\hbox{for all $\delta>0$, there is $C_\delta>0$ such that}\ w_0(x)\leq C_\delta e^{-(r_-(c)-\delta)x}\ \hbox{for}\ x>0,\\
&~~~~~~~~~~~~ w_0(x)\leq C_\delta e^{-(r_-(c)+\varepsilon)x}\ \hbox{for}\ x\in[(1-q\e)x_n,(1+q\e)x_n],
\end{array}
$$

(iii). and we have
$$
\begin{array}{rll}
\di{\int H(t,x-y)\biggl((u_0*\psi_0)(y)+C e^{C\vert x-y\vert/t}u_0(y)\biggl)\ dy}\leq &\di{\int H(t,x-y)w_0(y)\ dy}\\
\di{\int E(t,x,y)\biggl((u_0*\psi_0)(y)+C e^{C\vert x-y\vert/t}u_0(y)\biggl)\ dy}\leq &\di{\int H(t,x-y)w_0(y)\ dy}
\end{array}
$$
And thus, we start with
\begin{eqnarray}
\label{apr1904}
u_1(t,x)\leq \frac{Ce^t}{\sqrt{t}}\int_{\Rm } 
e^{-\frac{(x-y)^2}{4t}} w_0(y)\ dy. 
\end{eqnarray}
And, as in the proof of Lemma \ref{l3.1}, we are going to estimate $u_1(t,x)$ at the points
$$
t_n=\frac{x_n}{\sqrt{c^2-4}},\ \ \ z_n=ct_n.
$$
Observe that for $n$ sufficiently large, $L + M_0 < 0$, so $z_n > L+M_0$. Thus $z_n\notin\supp(a-1)$ and the estimate (\ref{apr1904}) applies. Let us decompose
\begin{eqnarray*}
&&u_1(t_n,z_n)=\frac{C e^{t_n}}{\sqrt{t_n}}\biggl(\int_{-\infty}^{0}+
\int_{0}^{(1-q\e)x_n}+\int^{(1+q\e)x_n}_{(1-q\e)x_n}+
\int_{(1+q\e)x_n}^{+\infty}\biggl) 
 {e^{-\frac{(z_n-y)^2}{4t_n}}} w_0(y)\ {dy} \\
&&~~~~~~~~~~~:=u_{11}(t_n,z_n)+u_{12}(t_n,z_n)+u_{13}(t_n,z_n)+u_{14}(t_n,z_n).
\end{eqnarray*}
As $ z_n-y \geq ct_n$ for $y\leq 0$, $t_n\ge 1$,
and $0\le w_0(y)\le 1$, we have
\begin{equation}\label{sep0306}
u_{11}(t_n,z_n)\leq Ce^{(1-\frac{c^2}4)t_n}\to 0\hbox{ as $n\to+\infty$},
\end{equation}
since $c>2$.
By Lemma \ref{l3.1} we have, for every $\delta>0$ 
\begin{equation}\label{apr1901}
u_{12}(t_n,z_n)\leq C_\delta\di{\int\limits_{0}^{(1-q\e)x_n}e^{t_n 
-\frac{(ct_n-y)^2}{4t_n}-(r_-(c)-\delta)y}\
\frac{dy}{\sqrt{t_n}}}.
\end{equation}
The  integrand above is maximized at 
the point 
\[
y_\delta=(c-2r+2\delta)t_n=
(\sqrt{c^2-4}+2\delta)t_n=x_n+\frac{2\delta}{\sqrt{c^2-4}}x_n,
\] 
that is $O(\delta x_n)$ close to $x_n$ -- this is, indeed, why $t_n$ was chosen
as above. Here we have used (\ref{apr2806}). As $y_\delta>x_n$, the integrand in (\ref{apr1901}) on the interval
$[0,(1-\eps q)x_n]$ is maximized at
 the upper limit, leading to
\begin{eqnarray*}
&&u_{12}(t_n,z_n)
\leq C\di{\int_{0}^{(1-q\e)x_n}e^{(1-(r_-(c)+ {q\e\sqrt{c^2-4}}/2)^2)t_n-(r_-(c)-\delta)(1-q\e)x_n}}
\farc{dy}{\sqrt{t_n}}\\
&&~~~~~~~~~~~~\leq
\di{{C}{\sqrt t_n}e^{[-q^2{(c^2-4)\e^2}/4+\delta(1-q\e)\sqrt{c^2-4}]t_n}   .          
}
\end{eqnarray*}
Recall that   $\e<1$. Hence, if we choose 
$\delta \leq\di\frac{q^2\e^2}{100}\sqrt{c^2-4}$ we have
$$
-q^2\frac{(c^2-4)\e^2}4+\delta(1-q\e)\sqrt{c^2-4}\leq-q^2\frac{(c^2-4)\e^2}8,
$$
and therefore
\beq
u_{12}(t_n,z_n)\leq 
C_\delta \sqrt{t_n}e^{-{q^2}\e^2(c^2-4)t_n/8}\to 0\hbox{ as $n\to +\infty$.}
\eeq
Consider now $u_{14}(t_n,z_n)$:
\begin{eqnarray*}
&&\!\!\!\!\!
u_{14}(t_n,z_n)\le \frac{Ce^{t_n}}{\sqrt{t_n}}\int\limits_{(1+q\e)x_n}^{+\infty} 
\!\!\!\!\! {e^{-\frac{(z_n-y)^2 }{4t_n}-(r_-(c)-\delta)y}}   {dy} 
=
Ce^{t_n}\left[\int_{(1+q\e)x_n}^{z_n}+\int_{z_n}^{+\infty}\right]
 \frac{e^{-\frac{|z_n-y|^2}{4t_n}-(r_-(c)-\delta)y}}
{\sqrt{ t_n}}  {dy}\\
&&~~~~~~~~~~~=u_{14}'(t_n,z_n)+u_{14}''(t_n,z_n).
\end{eqnarray*}
For $u_{14}''$ we have:
\begin{eqnarray*}
&&u_{14}''(t_n,z_n)=Ce^{t_n}\int_{z_n}^{+\infty} 
 \frac{e^{-\frac{(y-z_n)^2}{4t_n}-(r_-(c)-\delta)y}}
{\sqrt{ t_n}}  {dy}\le  Ce^{t_n-(r_-(c)-\delta)ct_n}=Ce^{-(r_-(c)^2-\delta)t_n}\to
0,
\end{eqnarray*}
as $n\to+\infty$,  while for $u_{14}'$ we have
\begin{eqnarray*}
&&u_{14}'(t_n,z_n)\le Ce^{t_n}\int_{(1+q\e)x_n}^{z_n}
\frac{e^{-\frac{(z_n-y)^2}{4t_n}-(r_-(c)-\delta)y}}{\sqrt{ t_n}}  {dy},
\end{eqnarray*}
and this term can be estimated exactly as $u_{12}(t_n,z_n)$.

We turn to 
$u_{13}(t_n,z_n)$ -- it is here that we use the crucial assumption (\ref{e3.3}). 
It follows from this bound on $w_0(y)$ inside the interval of integration
that
\begin{equation}\label{apr1902}
u_{13}(t_n,z_n)\leq C
\int\limits^{(1+q\e)x_n}_{(1-q\e)x_n}\frac{e^{t_n-\frac{(ct_n-y)^2-C|ct_n-y|}{4t_n} 
-(r_-(c)+\e/2)y}}{\sqrt{4\pi t_n}}\ dy
\le C
\int\limits^{(1+q\e)x_n}_{(1-q\e)x_n}\frac{e^{t_n-\frac{(ct_n-y)^2}{4t_n} 
-(r_-(c)+\e/2)y}}{\sqrt{4\pi t_n}}dy.
\end{equation}
Now, the maximum of the integrand is achieved at the point
\[
y_n=x_n-\frac{\eps}{\sqrt{c^2-4}}x_n.
\]
At the expense of possibly decreasing $q$
so that $q<\di 1/\sqrt{c^2-4}$, we have $y_n<(1-q\e)x_n$. Then the integrand in 
(\ref{apr1902}) is maximized at $y=(1-q\eps)x_n$,
and we have, for all $y\in[(1-q\e)x_n,(1+q\e)x_n]$:
\begin{eqnarray}
&&\di{-\frac{(ct_n-y)^2}{4t_n}-(r_-(c)+\frac\e2)y}\leq \di{-\frac{(ct_n-(1-q\e)x_n)^2}{4t_n}-(r_-(c)+\frac\e2)(1-q\e)x_n}\\
&&\leq\di{\biggl(-1-\frac\eps2\sqrt{c^2-4}+O(\e^2)\biggl)t_n}.\nonumber
\end{eqnarray}
This gives, for $\eps>0$ sufficiently small, 
\beq
u_{13}(t_n,z_n)\leq Cx_ne^{-\e t_n\sqrt{c^2-4}/4}
\eeq
and, all in all, we have the following upper bound for $u_1(t_n,z_n)$:
\beq
\label{e3.4}
u_1(t_n,z_n)\leq C \sqrt{t_n} e^{- \e t_n\sqrt{c^2-4}/4}+C_\delta \sqrt{t_n} e^{- 
{q^2}\e^2(c^2-4)t_n/8}.
\eeq

\subsubsection*{The estimate for $u_2(t_n,z_n)$}

The quantity $L+M_0$ is bounded from above by a universal constant,
so  
\begin{eqnarray}\label{oct404}
&&u_2(t_n,z_n)\le \farc{Ce^{t_n}}{\sqrt{t_n}}\int_{L-M_0}^\infty 
e^{-\farc{|z_n+y-2L|^2}{4t_n}}w_0(y)dy 
=  {Ce^{t_n}} \int_{(z_n-(L+M_0))/\sqrt{4t_n}}^\infty e^{-y^2}dy\leq Ce^{t_n-z_n^2/(4t_n)}
\nonumber\\
&&~~~~~~~~~~~~\leq Ce^{(1-c^2/4)t_n}.
\end{eqnarray}
This will decay exponentially fast since $c>2$.

\subsubsection*{Estimate of $u_3(t_n,z_n)$}

The last term we need to consider is the eigenvalue contribution:
\[
u_3(t,x)=e^{\blambda t}\phi_0(x)\int \phi_0(y)w_0(y)dy,
\]
and this is also easy: we have 
\begin{equation}\label{oct406}
u_3(t_n,z_n)\leq Ce^{\blambda t_n-\sqrt{\blambda-1}z_n}=Ce^{(\blambda -c\sqrt{\blambda-1})t_n},
\end{equation}
and this quantity will also decay exponentially fast because  $c>\blambda/\sqrt{\blambda-1}$.


\subsubsection* { The estimate for $u_4(t_n,z_n)$} 

We wish to show that $u_4(t_n,z_n)$ goes to $0$ as $n\to+\infty$  
slower than the first three terms. As the front is moving with the speed $c$, for any small $\delta>0$, 
there exists a large $x_\delta>0$ such that
$$
u(t,x)\geq\frac12\ \ \ \hbox{for $x\leq(c-\delta)t-x_\delta$ and $t\geq 0$}.
$$
By our assumption on $f(u)$ there is a constant $C > 0$ such that $u - f(u) \geq C$ for all $u \in [1/2,1]$. Therefore, as $a(x)\ge a_0>0$, we have 
\begin{eqnarray}
u_4(t_n,z_n) & \ge & a_0\int_{0}^{t_n}
\int_{\bbR}\frac{e^{t_n-s-\frac{(ct_n-y)^2}{4(t_n-s)}}}{\sqrt{4\pi(t_n-s)}}(u(s,y) - f(u(s,y))) 
dsdy \no \\
&\geq & C\int_{0}^{t_n}\int_{(c-\delta)s-x_\delta-1}^{(c-\delta)s-x_\delta}
\frac{e^{t_n-s- {(ct_n-y)^2}/{4(t_n-s)}}}{\sqrt{(t_n-s)}}\ dsdy.
\end{eqnarray}
The change of variables $y=(c-\delta)s-x_\delta+z$ in the last integral yields 
$$
u_4(t_n,z_n)\geq 
\frac{C}{\sqrt{t_n}}\int_{0}^{t_n}\int_{-1}^0
e^{t_n-s- {(c(t_n-s)+\delta s+x_\delta-z)^2}/{4(t_n-s)}}{dsdy} .
$$
We have, for $z\in(-1,0)$ and $0\le s<t_n-1$:
\begin{eqnarray*}
&&\Psi_\delta(s,t_n,z):=
\di{t_n-s-\frac{(c(t_n-s)+\delta s+x_\delta-z)^2}{4(t_n-s)}}\\
&&~~~~~~~~~~~~~~~=\di{(1-\frac{c^2}4)(t_n-s)-\farc{c\delta s}{2}
-\frac{\delta^2s^2}{4(t_n-s)}-2(x_\delta-z)\frac{c(t_n-s)+\delta s}{4(t_n-s)}
-\frac{(x_\delta-z)^2}{4(t_n-s)}}.
\end{eqnarray*}
We evaluate the integral on the time interval 
$(1-\gamma_1)t_n\le s \le (1-\gamma_2)t_n$ with 
$0<\gamma_2<\gamma_1\ll 1$ to be chosen. There is a constant 
$C_{\delta,\gamma}$ that depends on $\gamma_{1,2}$ and $\delta$ but not on 
$n$ such that for all $z\in[-1,0]$  and all
$s$ in this interval we have
$$
\begin{array}{rll}
\Psi_\delta(s,t_n,z)\geq&\di{(1-\frac{c^2}4)(t_n-s)-\farc{c\delta s}{2}
-\frac{\delta^2s^2}{4(t_n-s)}-C_{\delta,\gamma}}\\
\geq&\biggl((1-\di\frac{c^2}4)\gamma_1-\di\frac{c}2\delta-
\di\frac{\delta^2(1-\gamma_2)^2}{4\gamma_2}\biggl)t_n-C_{\delta,\gamma}
:=-A_{\delta,\gamma}t_n-C_{\delta,\gamma}.
\end{array}
$$
Therefore
\beq
\label{e3.6}
u_4(t_n,z_n)\geq {C}{\sqrt{t_n}}e^{-A_{\delta,\gamma}t_n-C_{\delta,\gamma}}.
\eeq
Gathering \eqref{e3.4}, \eqref{oct404}, \eqref{oct406} and \eqref{e3.6} we have, 
for a  constant $C>0$ depending only on $\delta$:
$$
u(t_n,z_n)\leq C_\delta\biggl(-e^{-A_{\delta,\gamma}t_n-C_{\delta,\gamma}}
+e^{-\e ct_n}+e^{-\eps t_n\sqrt{c^2-4}/4}+e^{(1-\frac{c^2}4+o(1))t_n}+
e^{-\frac{q^2}2\e^2(c^2-4)t_n}+e^{(\blambda-c\sqrt{\blambda-1})t_n}\biggl).
$$
Choosing  $\gamma_1$ and $\gamma_2$ small enough, and then $\delta=\gamma_2$   
makes the constant $A_{\delta,\gamma}$
arbitrarily small.   In particular, we may ensure that it is much smaller than 
the
coefficients in front of $t_n$ in the last five exponential terms above. This
 yields 
$$
u(t_n,z_n)<0
$$
for large $n$ which is the contradiction. \hfill $\qed$

\end{document}